\documentclass{article}
\usepackage{amssymb}
\usepackage{amsmath}

\newtheorem{theorem}{Theorem}

\newtheorem{corollary}[theorem]{Corollary}

\newtheorem{definition}[theorem]{Definition}

\newtheorem{lemma}[theorem]{Lemma}

\newtheorem{proposition}[theorem]{Proposition}
\newtheorem{remark}[theorem]{Remark}

\input{tcilatex}
\begin{document}

\begin{center}
{\LARGE A Pointwise a-priori Estimate\smallskip\ }

{\LARGE for the }$\overline{\partial }{\LARGE -}${\LARGE Neumann Problem
\smallskip \smallskip }

{\LARGE on Weakly Pseudoconvex Domains\footnote{\textbf{Mathematics Subject
Classification (2010):}
\par
Primary: 32W05, 35N15; Secondary: 32A26, 32T27.}}

\vspace{0.5in}

{\Large R. Michael Range\vspace{0.5in}}

\bigskip

\textbf{Abstract\bigskip }
\end{center}

The main result is a pointwise a-priori estimate for the $\overline{\partial 
}-$Neumann problem that holds on an \emph{arbitrary} weakly pseudoconvex
domain $D$. It is shown that for $(0,q)$ forms $f$ in the domain of the
adjoint $\overline{\partial }^{\ast }$ of $\overline{\partial }$, the
pointwise growth of the derivatives of each coefficient of $f$ with respect
to $\overline{z_{j}}$ and in complex tangential directions is carefully
controlled by the sum of the suprema of $f$, $\overline{\partial }f$, and $%
\overline{\partial }^{\ast }f$ over $D$. These estimates provide a pointwise
analogon of the classical basic estimate in the $L^{2}$ theory that has been
the starting point for all major work in this area involving $L^{2}$ and
Sobolev norm estimates for the complex Neumann and related operators. It is
expected that under suitable additional conditions on the boundary $bD$,
such as finite type, the estimates may be refined to yield H\"{o}lder
estimates. Such result have been known in case of \emph{strictly}
pseudoconvex domains for a long time, but well known obstructions---even in
case of finite type---have blocked progress on understanding pointwise
estimates for the $\overline{\partial }-$complex via techniques of integral
representations in sufficiently general cases. The obstructions are overcome
by introducing new Cauchy-Fantappi\'{e} kernels whose critical novel
ingredient is a non-holomorphic modification of the Levi polynomial of a
suitable defining function for the domain. The construction works on any
smoothly bounded pseudoconvex domain, and the resulting kernels reflect the
complex geometry of the boundary $bD$.

\begin{center}
\bigskip
\end{center}

\section{\protect\bigskip Introduction}

The $L^{2}$ theory of the $\overline{\partial }-$Neumann problem on
pseudoconvex domains has been highly developed for many years. In
particular, in 1979 J. J. Kohn introduced the technique of subelliptic
multipliers which led to the proof of subelliptic estimates in case the
boundary is of finite type ([1], [3], [8], [15]). The starting point for
these and other investigations has been the following basic estimate, valid
on any smoothly bounded pseudoconvex domain $D$. (See [6] or [8] for more
details.) Let us fix a point $P\in bD$ and a smooth orthonormal frame for $%
(1,0)$ forms $\omega _{1},$ $\omega _{2},...,\omega _{n}$ on a small
neighborhood $U$ of $P$ with $\omega _{n}=\gamma (\zeta )\partial r$, where $%
r$ is a defining function for $D$. Let $L_{1},...,L_{n}$ be the
corresponding dual frame for $(1,0)$ vector fields. One defines 
\begin{equation*}
\mathfrak{D}_{q}(D)=C_{(0,q)}^{\infty }(\overline{D})\cap dom(\overline{%
\partial }^{\ast }),
\end{equation*}%
and one denotes by $\mathfrak{D}_{qU}$ those forms in $\mathfrak{D}_{q}(D)$
which have compact support in $\overline{D}\cap U$. Then $f\in \mathfrak{D}%
_{qU}$ can be written as $\sum_{J}^{^{\prime }}f_{J}\overline{\omega }^{J}$,
where the summation is over strictly increasing $q$-tuples $J$. Since $f\in
dom(\overline{\partial }^{\ast }),$ one has $f_{J}=0$ on $bD\cap U$ whenever 
$n\in J$. In case $q=1,$ the \textquotedblleft $L^{2}$ basic
estimate\textquotedblright\ states that there exists a constant $C$ such
that 
\begin{equation*}
\sum_{j,~k}\left\Vert \overline{L_{j}}f_{k}\right\Vert ^{2}+\int_{bD\cap U}%
\mathcal{L(}r,\mathcal{\zeta };f^{\#})dS(\zeta )\leq C~\left[ \left\Vert 
\overline{\partial }f\right\Vert ^{2}+\left\Vert \overline{\partial }^{\ast
}f\right\Vert ^{2}+\left\Vert f\right\Vert ^{2}\right]
\end{equation*}%
for all $f=\sum f_{k}\overline{\omega _{k}}\in \mathfrak{D}_{1U}$, where $%
f^{\#}=(f_{1},...,f_{n})$. The norms here are the standard $L^{2}$ norms
over $D\cap U$ , and $\mathcal{L}$ is the Levi form of the defining function 
$r$ with respect to the frame $\{L_{1},...,L_{n}\}$. Since $f\in dom(%
\overline{\partial }^{\ast })$, one has $f_{n}=0$ on $bD$, so that
pseudoconvexity implies that $\mathcal{L(}r,\mathcal{\zeta };f^{\#})\geq 0$
on $bD$. Furthermore, it readily follows from $f_{n}\left\vert _{bD\cap
U}\right. =0$ that one then also has the estimate 
\begin{equation*}
\left\Vert f_{n}\right\Vert _{1}^{2}\leq C_{1}\left\Vert \overline{\partial }%
f_{n}\right\Vert ^{2}\leq C_{2}\left[ \left\Vert \overline{\partial }%
f\right\Vert ^{2}+\left\Vert \overline{\partial }^{\ast }f\right\Vert
^{2}+\left\Vert f\right\Vert ^{2}\right] \text{ .}
\end{equation*}%
Here $\left\Vert f_{n}\right\Vert _{1}$ is the full $1$-Sobolev norm, i.e., $%
\left\Vert f_{n}\right\Vert _{1}^{2}$ is the sum of the squares of the $%
L^{2} $ norms of all first order derivatives of $f_{n}$.\medskip

Over the years there has been much interest in obtaining corresponding
results involving pointwise and H\"{o}lder estimates. Techniques of integral
representations have been most successful on \emph{strictly} pseudoconvex
domains, where the Levi polynomial provides a simple explicit local
holomorphic support function. (See [12] for a systematic exposition.)
Holomorphic support functions also exist on convex domains, and some results
have been obtained in that setting in case of finite type ([2], [5]).
However, it has long been known that in general there are no analogous
holomorphic support functions, even in very simple pseudoconvex domains of
finite type ([9]). This obstruction has blocked any progress on these
questions in case of more general pseudoconvex domains.

Recently the author has introduced a \emph{non-holomorphic} modification of
the Levi polynomial to obtain new Cauchy-Fantappi\'{e} kernels on arbitrary
weakly pseudoconvex domains which reflect the complex geometry of the
boundary and which satisfy some significant partial estimates ([14]). In
this paper we use the new kernels in the integral representation formula
developed by I. Lieb and the author in the strictly pseudoconvex case (see
[10]) to prove a pointwise analogon of the classical basic $L^{2}$ estimate,
as follows. This result was already announced in [13]. We define 
\begin{equation*}
\mathfrak{D}_{q}^{k}(D)=C_{(0,q)}^{k}(\overline{D})\cap dom(\overline{%
\partial }^{\ast })
\end{equation*}%
for $k=1,2,...$, and we denote by $\mathfrak{D}_{qU}^{k}$ those forms in $%
\mathfrak{D}_{q}^{k}(D)$ that have compact support in $\overline{D}\cap U.$
We shall use the frames $\omega _{1},$ $\omega _{2},...,\omega _{n}$ and $%
L_{1},...,L_{n}$ as above. Vectorfields $V$ act on forms coefficientwise,
i.e., if $f=\sum_{J}^{^{\prime }}f_{J}\overline{\omega }^{J},$ then $%
V(f)=\sum_{J}^{^{\prime }}V(f_{J})\overline{\omega }^{J}$. For a $C^{1}$
form $f$ of type $(0,q)$ on $\overline{D}$ we define the norm 
\begin{equation*}
Q_{0}(f)=\left\vert f\right\vert _{0}+\left\vert \overline{\partial }%
f\right\vert _{0}+\left\vert \vartheta f\right\vert _{0},
\end{equation*}%
where $\vartheta $ is the formal adjoint of $\overline{\partial }$, and $%
\left\vert \varphi \right\vert _{0}$ denotes the sum of the supremum norms
over $D$ of the coefficients of $\varphi .$ For $0<\delta \leq 1,$ $%
\left\vert \varphi \right\vert _{\delta }$ denotes the corresponding H\"{o}%
lder norm of order $\delta ,$

\textbf{Main Theorem. }\emph{There exists an integral operator }$%
S^{bD}:C_{(0,q)}(bD)\rightarrow C_{(0,q)}^{\infty }(D)$ \emph{which has the
following properties. If }$bD$\emph{\ is (Levi) pseudoconvex in a
neighborhood }$U$\emph{\ of the point} $P\in bD$ \emph{and} \emph{if }$U$%
\emph{\ is sufficiently small, there exist constants }$C_{\delta }$ \emph{%
depending on} $\delta >0$, \emph{so that one has the following uniform
estimates for all }$f\in \mathcal{D}_{qU}^{1}$\emph{, }$1\leq q\leq n$, 
\emph{and for }$z\in D\cap U$\emph{\ :}

\emph{1) }$\left\vert f-S^{bD}(f)\right\vert _{\delta }\leq C_{\delta }\cdot
Q_{0}(f)$ \emph{for any }$\delta <1.$

\emph{2) }$\left\vert \overline{L_{j}}S^{bD}(f)(z)\right\vert \leq C_{\delta
}\cdot dist(z,bD)^{\delta -1}\cdot Q_{0}(f)$\emph{\ for }$j=1,...,n$\emph{\
and any }$\delta <1/2$\emph{;}

\emph{3) }$\left\vert L_{j}S^{bD}(f)(z)\right\vert \leq C_{\delta }\cdot
dist(z,bD)^{\delta -1}\cdot Q_{0}(f)$\emph{\ for }$j=1,...,n-1$\emph{\ and
any }$\delta <1/3$\emph{;}

\emph{\ Furthermore, if }$f_{J}\overline{\omega }^{J}$ \emph{is a normal
component of }$f$\emph{\ with respect to the frame }$\overline{\omega }%
_{1},...,\overline{\omega }_{n}$,\emph{\ one has}

\begin{equation*}
\left\vert f_{J}\right\vert _{\delta }\leq C_{\delta }Q_{0}(f)\text{\emph{\
for any }}\delta <1/2\text{ \emph{if }}n\in J\text{.}
\end{equation*}

Note that if one had an estimate analogous to 3) also for the normal
derivative $L_{n}S^{bD}(f)(z)$ for \emph{some} $\delta >0$ (with $\delta
<1/3 $), standard results would imply the H\"{o}lder estimate $\left\vert
S^{bD}(f)\right\vert _{\delta }\leq C_{\delta }Q_{0}(f)$; by using 1) one
therefore would obtain an estimate 
\begin{equation*}
\left\vert \,f\,\right\vert _{\delta }\leq C_{\delta }Q_{0}(f)\text{\emph{,}}
\end{equation*}%
i.e., the H\"{o}lder analogon of a subelliptic estimate. It is known that
such an estimate does not hold on arbitrary pseudoconvex domains. On the
other hand, the Main Theorem provides a starting point in a general setting
which, combined with additional suitable properties of the boundary such as
finite type, might be useful to obtain appropriate estimates for $%
L_{n}S^{bD}(f).$ In particular, the author is investigating analogs of
Kohn's subelliptic multipliers in the integral representation setting
underlying the Main Theorem. (See [13] for an outline of such potential
applications.)

\section{Integral Representations}

We shall briefly recall some fundamentals of the integral representation
machinery. We shall follow the terminology and notations from [12], where
full details may be found. A (kernel) generating form $W(\zeta ,z)$ for the
smoothly bounded domain $D\subset \mathbb{C}^{n}$ is a $(1,0)$ form $%
W=\sum_{j=1}^{n}w_{j}d\zeta _{j}$ defined on $bD\times D$ with coefficients
of class $C^{1}$ which satisfies $\sum w_{j}(\zeta _{j}-z_{j})=1$. For $%
0\leq q\leq n-1,$ the associated Cauchy-Fantappi\'{e} (= CF) form of order $%
(0,q)$ is defined by 
\begin{equation*}
\Omega _{q}(W)=c_{nq}W\wedge (\overline{\partial _{\zeta }}W)^{n-q-1}\wedge (%
\overline{\partial _{z}}W)^{q}\;\footnote{$c_{nq}=\frac{(-1)^{q(q-1)/2}}{%
(2\pi i)^{n}}\left( 
\begin{array}{c}
n-1 \\ 
q%
\end{array}%
\right) $}\text{.}
\end{equation*}%
$\Omega _{q}(W)$ is a double form on $bD\times D$ of type $(n,n-q-1)$ in $%
\zeta $ and type $(0,q)$ in $z$. One also sets $\Omega _{-1}(W)=\Omega
_{n}(W)=0$.

With $\beta =\left\vert \zeta -z\right\vert ^{2}$, the form $B=\partial
\beta /\beta =\sum_{j=1}^{n}\overline{(\zeta _{j}-z_{j})}/\left\vert \zeta
-z\right\vert ^{2}d\zeta _{j}$ is the generating form for the
Bochner-Martinelli-Koppelman (= BMK) kernels. One has the following BMK
formula. (Here and in the following the integration variable is always $%
\zeta $.)

\emph{If }$f\in C_{(0,q)}^{1}(\overline{D}),$\emph{\ then} 
\begin{equation}
f(z)=\int_{bD}f(\zeta )\wedge \Omega _{q}(B)-\overline{\partial _{z}}%
\int_{D}f(\zeta )\wedge \Omega _{q-1}(B)-\int_{D}\overline{\partial _{\zeta }%
}f(\zeta )\wedge \Omega _{q}(B)\text{ \ \emph{for }}z\in D\text{.}
\label{BMK}
\end{equation}

The next formula, due to W. Koppelman, describes how to replace $\Omega
_{q}(B)$ by some other CF kernel $\Omega _{q}(W)$ on the boundary $bD.$
Since $\Omega _{n}(B)\equiv 0$, we shall assume $q<n$ from now on. Given any
generating form $W$ on $bD\times D$, one has

\begin{equation*}
f(z)=\int_{bD}f(\zeta )\wedge \Omega _{q}(W)+\overline{\partial _{z}}%
~T_{q}^{W}(f)+T_{q+1}^{W}(\overline{\partial }f)\text{ \emph{for} }f\in
C_{(0,q)}^{1}(\overline{D})\text{ \emph{and} }z\in D\text{.}
\end{equation*}%
Here the integral operator $T_{q}^{W}:C_{(0.q)}^{1}(\overline{D})\rightarrow
C_{(0,q-1)}(D)$ is defined by 
\begin{equation*}
T_{q}^{W}(f)=\int_{bD}f\wedge \Omega _{q-1}(W,B)-\int_{D}f(\zeta )\wedge
\Omega _{q-1}(B)
\end{equation*}%
for any $0\leq q<n$, where the \textquotedblleft
transition\textquotedblright\ kernels $\Omega _{q-1}(W,B)$ involve explicit
expressions in terms of $W$ and $B$ which will be recalled later on.

\textbf{Remark.} In case $D$ is strictly pseudoconvex, Henkin and Ramirez
constructed a generating form $W^{HR}(\zeta ,z)$ for $D$ which is
holomorphic in $z$, so that $\Omega _{q}(W^{HR})=0$ on $bD$ for $q\geq 1$.
Consequently, if $f$ is a $\overline{\partial }$-closed $(0,q)$ form on $%
\overline{D}$, one has $f=\overline{\partial _{z}}~T_{q}^{W^{HR}}(f)$, with
an explicit solution operator $T_{q}^{W^{HR}}$. Based on the critical
information that the Levi form of the boundary is positive definite in this
case, it is well-known that this solution operator is bounded from $%
L^{\infty }$ into $\Lambda _{1/2}$. Furthermore, one also has the a-priori H%
\"{o}lder estimate $\left\vert f\right\vert _{1/2}\leq CQ_{0}(f)$\emph{\ }%
for all $f\in \mathcal{D}_{qU}^{1}$. (See [11]). Attempts to prove
corresponding estimates on more general domains ultimately run into the
obstruction of the Kohn - Nirenberg example mentioned above ([9]), i.e., in
general it is not possible to find a corresponding reasonably explicit \emph{%
holomorphic} generating form on weakly pseudoconvex domains---even if of
finite type---except under very restrictive geometric conditions.

In case $f\in \mathfrak{D}_{q}^{1}(D),$ one may transform formula (\ref{BMK}%
) into

\begin{equation*}
f=\int_{bD}f\wedge \Omega _{q}(B)+(\overline{\partial }f,\overline{\partial }%
\omega _{q})+(\vartheta f,\vartheta \omega _{q-1})\text{,}
\end{equation*}%
where $\omega _{q}$ denotes the fundamental solution of $\square $ on $(0,q)$
forms, $\vartheta $ denotes the formal adjoint of $\overline{\partial }$, so
that $\vartheta f=\overline{\partial }^{\ast }f$, and $(\bullet ,\bullet )$
denotes the standard $L^{2}$ inner product of forms over $D$. (See [10,
12]). The fundamental solution $\omega _{q}$ is an isotropic kernel whose
regularity properties are well understood. In particular, the operator 
\begin{equation*}
S^{iso}:f\rightarrow S^{iso}(f)=(\overline{\partial }f,\overline{\partial }%
\omega _{q})+(\vartheta f,\vartheta \omega _{q-1})
\end{equation*}
satisfies a H\"{o}lder estimate%
\begin{equation}
\left\vert S^{iso}(f)\right\vert _{\delta }\leq C_{\delta }Q_{0}(f)\text{
for all }f\in C_{(0,q)}^{1}(D)\text{ and any }\delta <1.  \label{iso}
\end{equation}%
Consequently, the essential information regarding all pointwise estimations
is contained in the boundary integral $S^{bD}(f)=\int_{bD}f\wedge \Omega
_{q}(B)$. Note that the kernel of $\Omega _{q}(B)$ is isotropic; it treats
derivatives in all directions equally, and direct differentiation under the
integral in $\int_{bD}f\wedge \Omega _{q}(B)$ leads to an expression that
will in general blow up like $dist(z,bD)^{-1}.$ So this general
representation of the operator $S^{bD}$ does not provide any useful
information.

Note that since $\Omega _{n}(B)\equiv 0,$ the Main Theorem trivially holds
with $S^{bD}\equiv 0$ in case $q=n$.

By the Koppelman formulas, given any generating form $W$ on $bD\times D$,
one can transform $S^{bD}(f)$ into%
\begin{equation}
S^{bD}(f)=\int_{bD}f\wedge \Omega _{q}(W)+\int_{bD}\overline{\partial }%
f\wedge \Omega _{q}(W,B)+\int_{bD}f\wedge \overline{\partial }_{z}\Omega
_{q-1}(W,B).  \label{Kop}
\end{equation}

The proof of the Main Theorem relies on formula (\ref{Kop}) on a weakly
pseudoconvex domain $D\subset \subset \mathbb{C}^{n}$, applied to the \emph{%
non-holomorphic} generating form $W^{\mathcal{L}}(\zeta ,z)$ introduced in
[14]. Let as briefly recall the key properties of $W^{\mathcal{L}}(\zeta ,z)$%
. Given a sufficiently small neighborhood $U=U(P)$, on $(bD\cap U)\times
(D\cap U)$ the form \medskip $W^{\mathcal{L}}(\zeta ,z)$ is represented
explicitly by

\begin{equation*}
W^{\mathcal{L}}(\zeta ,z)=\frac{\sum_{j=1}^{n}g_{j}(\zeta ,z)d\zeta _{j}}{%
\Phi _{K}(\zeta ,z)}\text{,}
\end{equation*}%
where $\Phi _{K}(\zeta ,z)=\sum_{j=1}^{n}g_{j}(\zeta ,z)(\zeta _{j}-z_{j})$
for $\zeta \in bD$. The (non-holomorphic) support function $\Phi _{K}$ is
defined by 
\begin{equation*}
\Phi _{K}(\zeta ,z)=F^{(r)}(\zeta ,z)-r(\zeta )+K\left\vert \zeta
-z\right\vert ^{3}\text{,}
\end{equation*}%
where $F^{(r)}(\zeta ,z)$ is the Levi polynomial of a suitable defining
function $r$, and $K>0$ is a suitably chosen large constant. We note that $%
W^{\mathcal{L}}(\zeta ,z)$ is $C^{\infty }$ in $z$ for $z\neq \zeta $.
Recall from [14] that the neighborhood $U$, the constant $K$, and $%
\varepsilon >0$ can be chosen so that for all $\zeta ,z$ $\in \overline{D}%
\cap U$ with $\left\vert \zeta -z\right\vert <\varepsilon $ one has

\begin{eqnarray}
\left\vert \Phi _{K}(\zeta ,z)\right\vert &\gtrsim &\left[ \left\vert \func{%
Im}F^{(r)}(\zeta ,z)\right\vert +\left\vert r(\zeta )\right\vert +\left\vert
r(z)\right\vert +\right.  \label{estimate} \\
&&\left. +\mathcal{L}(r,\zeta ;\pi _{\zeta }^{t}(\zeta -z))+K\left\vert
\zeta -z\right\vert ^{3}\right] \text{. }  \notag
\end{eqnarray}

Here $\pi _{\zeta }^{t}(\zeta -z)$ denotes the projection of $(\zeta -z)$
onto the complex tangent space of the level surface $M_{r(\zeta )}$ through
the point $\zeta $, and $\mathcal{L}(r,\zeta ;\pi _{\zeta }^{t}(\zeta -z))$
denotes the Levi form of $r$ at the point $\zeta $. As shown in [14], the
defining function $r$ can be chosen so that $\mathcal{L}(r,\zeta ;\pi
_{\zeta }^{t}(\zeta -z))\geq 0$ for all $\zeta $ $\in \overline{D}\cap U$.
We also recall that---as in the classical strictly pseudoconvex case---$%
r(\zeta )$ and $\func{Im}F^{(r)}(\zeta ,z)$ can be used as (real)
coordinates in a neighborhood of a fixed point $z.$

Note that since $\left\vert \zeta -z\right\vert ^{3}$ is real and symmetric
in $\zeta $ and $z$, it follows from the known case $K=0$ (see [12], for
example) that if one defines $\Phi _{K}^{\ast }(\zeta ,z)=\overline{\Phi
_{K}(z,\zeta )}$, one has the approximate symmetry 
\begin{equation}
\Phi _{K}^{\ast }-\Phi _{K}=\mathcal{E}_{3}\footnote{$\mathcal{E}_{j}$
denotes a smooth expression which satisfies $\left\vert \mathcal{E}%
_{j}\right\vert \lesssim \left\vert \zeta -z\right\vert ^{j}$.}.  \label{sym}
\end{equation}

In the following we shall simplify notation by dropping the subscript $K$,
i.e., we will write $\Phi $ instead of $\Phi _{K}$

\bigskip For $0\leq q<n$ we shall thus consider the integral representation
formula%
\begin{equation}
f=S^{bD}(f)+S^{iso}(f)\text{ for }f\in \mathfrak{D}_{0,q}^{1}(\overline{D})
\end{equation}%
where the boundary operator $S^{bD}$ is given by%
\begin{equation}
S^{bD}(f)=\int_{bD}f\wedge \Omega _{q}(W^{\mathcal{L}})+\int_{bD}\overline{%
\partial }f\wedge \Omega _{q}(W^{\mathcal{L}},B)+\int_{bD}f\wedge \overline{%
\partial }_{z}\Omega _{q-1}(W^{\mathcal{L}},B)\text{.}  \label{bDLevi}
\end{equation}%
for $f\in \mathfrak{D}_{0,q}^{1}(\overline{D})$. Corresponding formulas hold
locally on $U\cap bD$ whenever the boundary is Levi pseudoconvex in $U$. It
is then clear that property 1) in the Main Theorem is satisfied. The main
difficulty involves establishing the estimates 2) and 3).

The proof of the Main Theorem involves a careful analysis of the boundary
integrals in formula (\ref{bDLevi}). In contrast to [10], for the most part
we shall deal directly with the integrals over $bD$, thereby simplifying the
analysis. However, for the most critical terms we will need to apply Stokes'
theorem and introduce the Hodge $\ast $ operator as in [10] to transform the
integrals into standard $L^{2}$ inner products of forms over $D\cap U$, and
exploit certain approximate symmetries in the kernels.

\section{\protect\bigskip The Integral $\protect\int_{bD}f\wedge \Omega
_{q}(W^{\mathcal{L}})$}

In case $D$ is strictly pseudoconvex, $W^{\mathcal{L}}$ can be chosen to be
holomorphic in $z$ for $\zeta $ close to $z$, so that the estimations become
trivial if $q\geq 1$, since then $\overline{\partial }_{z}W=0$ near the
singularity. In the general case considered here, this integral needs to be
carefully estimated as well. The analysis of this integral involves
straightforward modifications of the case $q=0$ discussed in [14], as
follows.

We only consider $\zeta ,z$ with $\left\vert \zeta -z\right\vert
<\varepsilon /2$, so that we can use the explicit form of $W^{\mathcal{L}%
}=g/\Phi $ recalled above, and the local frames $\{\omega _{1},...,\omega
_{n}\}$ and $\{L_{1},...,L_{n}\}$. Recall that for $j=0,1,2$ an expression $%
\mathcal{E}_{j}^{\#}$ denotes a form which is smooth for $\zeta \neq z$ and
that satisfies a uniform estimate $\left\vert \mathcal{E}_{j}^{\#}\right%
\vert \lesssim \left\vert \zeta -z\right\vert ^{j},$ and whose precise
formula may change from place to place. While $\Phi $ is not holomorphic in $%
z,$ one has $\overline{\partial }_{z}\Phi =\mathcal{E}_{2}^{\#}$ ;
furthermore, one has $L_{j}^{z}\Phi =\mathcal{E}_{1}^{\#}$ for $j<n$, while $%
L_{n}^{z}\Phi \neq 0$ at $\zeta =z$.

By the properties of CF forms, on $bD$ one has 
\begin{equation*}
\Omega _{q}(W^{\mathcal{L}})=c_{nq}\frac{g\wedge (\overline{\partial }%
_{\zeta }g)^{n-q-1}\wedge (\overline{\partial }_{z}g)^{q}}{\Phi ^{n}}
\end{equation*}

The coefficients $g_{j}$ of $g=\sum g_{j}d\zeta _{j}$ are given by 
\begin{equation*}
g_{j}=\partial r/\partial \zeta _{j}-1/2\sum_{k}\partial ^{2}r/\partial
\zeta _{j}\partial \zeta _{k}(\zeta _{k}-z_{k})+\mathcal{E}_{2}^{\#}\text{.}
\end{equation*}%
The form of $g$ implies that 
\begin{eqnarray*}
g &=&\partial r(\zeta )+\mathcal{E}_{1}+\mathcal{E}_{2}^{\#}\text{, } \\
\overline{\partial }_{\zeta }g &=&\overline{\partial }\partial r(\zeta )+%
\mathcal{E}_{1}^{\#}\text{, and} \\
\overline{\partial }_{z}g &=&\mathcal{E}_{1}^{\#}\text{.}
\end{eqnarray*}

It follows readily that for $0\leq t\leq n-1$ one has 
\begin{equation*}
g\wedge (\overline{\partial _{\zeta }}g)^{t}=\partial r(\zeta )\wedge
\sum_{k=0}^{t}[\overline{\partial }\partial r(\zeta )]^{k}(\mathcal{E}%
_{1}^{\#})^{t-k}+\sum_{k=0}^{t}[\overline{\partial }\partial r(\zeta )]^{k}(%
\mathcal{E}_{1}^{\#})^{t-k+1}\text{,}
\end{equation*}%
where $(\mathcal{E}_{1}^{\#})^{s}$ denotes a generic form of appropriate
degree whose coefficients are products of $s$ terms of type $\mathcal{E}%
_{1}^{\#}$ in case $s\geq 1$, or a term of type $\mathcal{E}_{0}^{\#}$ for $%
s=0,-1$.

Note that since $\iota ^{\ast }(\omega _{n}\wedge \overline{\omega _{n}})=$ $%
0$ on $bD$, the pull-back of $\partial r(\zeta )\wedge \lbrack \overline{%
\partial }\partial r(\zeta )]^{k}$ to $bD$ involves only \emph{tangential}
components $tan[\overline{\partial }\partial r(\zeta )]$ of $\overline{%
\partial }\partial r(\zeta )$, while the pull-back of $[\overline{\partial }%
\partial r(\zeta )]^{k}$ alone will involve exterior products of at least $%
k-1$ different tangential components.

When estimating integrals involving these expressions we shall make use of
the fact that---in suitable $z-$diagonalizing coordinates (see [14])---each 
\emph{tangential }component $tan[\overline{\partial }\partial r(\zeta )]$ of 
$\overline{\partial }\partial r(\zeta )$ in the numerator of the kernel
reduces the order of vanishing of the corresponding factor $\Phi $ in the
denominator from $3$ to an estimate $\gtrsim \left\vert \zeta
_{l}-z_{l}\right\vert ^{2}$, i.e., 
\begin{equation}
\left\vert tan[\overline{\partial }\partial r(\zeta )]/\Phi \right\vert
\lesssim 1/(\left\vert r\left( z\right) \right\vert +\left\vert \zeta
_{l}-z_{l}\right\vert ^{2}),  \label{tanest}
\end{equation}%
where $\zeta _{l}$ is an appropriate complex tangential coordinate.
Similarly, 
\begin{equation*}
\left\vert \mathcal{E}_{1}^{\#}/\Phi \right\vert \lesssim 1/(\left\vert
r\left( z\right) \right\vert +\left\vert \zeta -z\right\vert ^{2}).
\end{equation*}%
In order to keep track of these estimates, we introduce forms \emph{$%
\mathcal{L[\mu ]}$} of \emph{Levi weight $\mathcal{\mu }$ }as follows. If $%
\mu \geq 1,$ we say \emph{$\mathcal{L[\mu ]}$} has Levi weight $\mu $ if
each summand of \emph{$\mathcal{L[\mu ]}$} contains at least $\mu $ factors
which either are (different) purely tangential components of $\overline{%
\partial }\partial r(\zeta )$, or of type $\mathcal{E}_{1}^{\#}.$ We also
set \emph{$\mathcal{L[\mu ]=E}_{0}$ }if $\mu \leq 0$. It then follows that 
\begin{equation*}
g\wedge (\overline{\partial _{\zeta }}g)^{t}=\mathcal{L}[t]
\end{equation*}%
on the boundary, and consequently the numerator of $\Omega _{q}(W^{\mathcal{L%
}})$ satisfies 
\begin{equation}
g\wedge (\overline{\partial }_{\zeta }g)^{n-q-1}\wedge (\overline{\partial }%
_{z}g)^{q}=\mathcal{L}[n-1]\text{. }  \label{num}
\end{equation}

\begin{proposition}
For any $q$ with $0\leq q\leq n-1$ the operator $T_{q}^{\mathcal{L}%
}:C_{(0,q)}(bD)\rightarrow C_{(0,q)}^{\infty }(D)$ defined by 
\begin{equation*}
T_{q}^{\mathcal{L}}f(z)=\int_{bD}f(\zeta )\wedge \Omega _{q}(W^{\mathcal{L}%
})(\zeta ,z)
\end{equation*}%
satisfies the estimates%
\begin{equation}
\left\vert \overline{L_{j}^{z}}(T_{q}^{\mathcal{L}}f(z))\right\vert \leq
C_{\delta }\left\vert f\right\vert _{0}dist(z,bD)^{\delta -1}\text{ for }%
\delta <2/3\text{ and }1\leq j\leq n\text{ }  \label{z-bar est}
\end{equation}%
and 
\begin{equation}
\left\vert L_{j}^{z}(T_{q}^{\mathcal{L}}f(z))\right\vert \leq C_{\delta
}\left\vert f\right\vert _{0}dist(z,bD)^{\delta -1}\text{ for }\delta <1/3%
\text{ and }j\leq n-1\text{ }  \label{tan est}
\end{equation}%
for suitable constants $C_{\delta }.$
\end{proposition}

Given the estimation (\ref{num}) of the numerator of $\Omega _{q}(W^{%
\mathcal{L}})$, the proof given in [14] for the case $q=0$ and for the
derivatives $\overline{L_{j}^{z}}$ carries over to the general case. To
prove the estimate (\ref{tan est}), one uses $L_{j}^{z}\Phi =\mathcal{E}%
_{1}^{\#}$ for $j\leq n-1$, which implies that $\left\vert L_{j}^{z}\Phi
/\Phi \right\vert \preceq dist(z,bD)^{-2/3}$. The estimations then proceed
as in [14].

\textbf{Remark.} There is no corresponding estimate for the differentiation
with respect to $L_{n}^{z}$, i.e., in the normal direction, since $%
L_{n}^{z}\Phi \neq 0$ for $\zeta =z$, so the operator $T_{q}^{\mathcal{L}}$
is not smoothing, i.e., there is no H\"{o}lder estimate 
\begin{equation*}
\left\vert T_{q}^{\mathcal{L}}f\right\vert _{\delta }\leq C_{\delta
}\left\vert f\right\vert _{0}\text{ for any }\delta >0\text{.}
\end{equation*}

Proposition 1 provides a \emph{partial} smoothing property, as follows.

\begin{definition}
We say that a kernel $\;\Gamma (\zeta ,z)$, or the integral operator $%
T_{\Gamma }:C_{\ast }(\overline{D})\rightarrow C_{\ast }^{1}(D)$ defined by
it, is $\overline{z}-$smoothing of order $\delta >0$ if $T_{\Gamma }$
satisfies the estimates (\ref{z-bar est}). Similarly, we say that $\Gamma $
(or $T_{\Gamma }$) is tangentially smoothing of order $\delta $ if $%
T_{\Gamma }$ satisfies the estimates (\ref{tan est}) for $L_{j}^{z}$ and $%
\overline{L_{j}^{z}}$ for $j=1,...,n-1$.
\end{definition}

Here $C_{\ast }$ denotes spaces of forms of appropriate type.

\section{\protect\bigskip Boundary Admissible Kernels}

Before proceeding with the analysis of the integrals involving the
transition kernels, we introduce admissible kernels and their weighted order
by suitably modifying corresponding notions from [10]. We say that a kernel $%
\Gamma (\zeta ,z)$ defined on $bD\times \overline{D}-\{(\zeta ,\zeta ):\zeta
\in bD\}$ is simple \emph{admissible }if for each $P\in bD$ there exists a
neighborhood $U$ of $P$, such that on $(bD\cap U)\times (\overline{D}\cap U)$
there is a representation of the form%
\begin{equation*}
\Gamma =\frac{\mathcal{L}[\mu ](\mathcal{E}_{1}^{\#})^{j}}{\Phi ^{t_{1}}%
\mathcal{\beta }^{t_{0}}}\text{,}
\end{equation*}%
where all exponents are $\geq 0.$ $j$ may be non-integer, in which case $(%
\mathcal{E}_{1}^{\#})^{j}$ denotes a form which is estimated by $C\left\vert
\zeta -z\right\vert ^{j}$. Such a representation is said to have (weighted) 
\emph{boundary} order\emph{\ } $\geq \lambda $ ($\lambda \in \mathbb{R}$)
provided\medskip

i) if $t_{1}\geq 1$, then

\begin{equation*}
2n-1+j-1-2\max (0,\min (t_{1}-1,\mu ))-3\max (t_{1}-1-\mu ,0)-2t_{0}\geq
\lambda
\end{equation*}%
if $\mu \geq 1,$ and%
\begin{equation*}
2n-1+j-1-3\max (t_{1}-1,0)-2t_{0}\geq \lambda \text{ }
\end{equation*}%
if $\mu \leq 0$,\medskip

or ii) if $t_{1}=0$, 
\begin{equation*}
2n-1+j-2t_{0}\geq \lambda \text{.}
\end{equation*}%
\medskip This definition of order takes into account that the dimension of $%
bD$ is $2n-1,$ and that \emph{one} factor $\Phi $ may be counted with weight
1, since by estimate (\ref{estimate}) one has $\left\vert \Phi \right\vert
\gtrsim \left\vert \func{Im}F^{r}\right\vert $, and $\func{Im}F^{r}(\cdot
,z) $ serves as a local coordinate on the boundary in a neighborhood of $z$.

A kernel $\Gamma $ is admissible of boundary order $\geq \lambda $ if it is
a finite sum of simple admissible kernels with representations of boundary
order $\geq \lambda $.

The results in the previous section show that $\Omega _{q}(W^{\mathcal{L}})$
is admissible of boundary order $\geq 0.$

As in the strictly pseudoconvex case considered in [10], an admissible
kernel $\Gamma $ of boundary order $\lambda \geq 1$ is smoothing of some
positive order $\delta $. This will follow from an estimate 
\begin{equation*}
\left\vert V^{z}T_{\Gamma }(f)(z)\right\vert \leq C_{\delta }\left\vert
f\right\vert _{0}dist(z,bD)^{-1+\delta }.
\end{equation*}%
for any vector field $V^{z}$ of unit length acting in $z$. On the other
hand, admissible kernels of boundary order $\lambda =0$ are not smoothing in
general.

More precisely, we have the following result.

\begin{theorem}
Let $\Gamma _{\lambda }$ be an admissible kernel of boundary order $\geq
\lambda $, and let 
\begin{equation*}
J_{\lambda }(z)=\int_{bD}\left\vert \Gamma _{\lambda }(\zeta ,z)\right\vert
dS(\zeta )\text{.}
\end{equation*}%
$J_{\lambda }(z)$ has the following properties.
\end{theorem}

\begin{enumerate}
\item[a)] \textit{If }$\lambda >0$\textit{, then }$\sup_{D}J_{\lambda
}(z)<\infty $\textit{.}

\item[b)] If $\lambda =0$, then $J_{0}(z)\lesssim dist(z,bD)^{-\alpha }$ for
any $\alpha >0$.

\item[c)] \textit{If }$\lambda \geq 1$\textit{, then }$\Gamma _{\lambda }$%
\textit{\ is smoothing of order }$\delta $\textit{\ for any }$\delta <1/3$,
tangentially smoothing of order $\delta <2/3$, and $\overline{z}$-smoothing
of order $\delta <1.$

\item[d)] \textit{If }$\lambda \geq 2$\textit{, then }$\Gamma _{\lambda }$%
\textit{\ is smoothing of order }$\delta $\textit{\ for any }$\delta <2/3.$
\end{enumerate}

\textit{Proof.} a) was essentially proved in [14] for the kernel $\Omega
_{0}(W^{\mathcal{L}})$. The general case follows by the same arguments. Part
b) follows from a) by noting that $dist(z,bD)^{\alpha }\leq \left\vert \zeta
-z\right\vert ^{\alpha }$ for $\zeta \in bD$. For c) note that given a
vector field $V^{z}$, all terms in $V^{z}\Gamma _{\lambda }$ are of boundary
order $\geq \lambda -1\geq 0$, except those where differentiation is applied
to $\Phi $; in that case use $V^{z}(\Phi ^{-s})=(\Phi ^{-s})[\mathcal{E}%
_{0}^{\#}/\Phi ]$ and $\left\vert 1/\Phi \right\vert \lesssim \left\vert
r(z)\right\vert ^{-2/3}1/\left\vert \zeta -z)\right\vert $ to see that $%
V^{z}\Gamma $ is estimated by $\left\vert r(z)\right\vert ^{-2/3}$
multiplied with a kernel of order $\geq 0.$ Similarly, in case $V^{z}$ is
tangential, one can replace $\left\vert r(z)\right\vert ^{-2/3}$ by $%
\left\vert r(z)\right\vert ^{-1/3}$, and in case $V^{z}=\overline{L_{j}^{z}}$%
, one uses $\overline{L_{j}^{z}}(\Phi ^{-s})=(\Phi ^{-s})[\mathcal{E}%
_{2}^{\#}/\Phi ]$ to see that $\overline{L_{j}^{z}}\Gamma _{\lambda }$ is of
order $\geq 0.$ The required estimates then follow from b). Finally, d)
follows by appropriately modifying the proof of c).$\blacksquare $

Most significant for this paper is the analysis of kernels of order $0$.
Such kernels are not smoothing in general. However, as we saw for $\Omega
_{q}(W^{\mathcal{L}})$, it turns out that in many cases they are at least $%
\overline{z}$-smoothing and tangentially smoothing of some positive order.
On the other hand, one readily checks that kernels of type such as $\mathcal{%
E}_{1}^{\#}/\beta ^{n}$ (e.g., those appearing in the BMK kernels) or $%
1/(\Phi \beta ^{n-1})$, which are of boundary order $0$, do not give
preference to tangential or $\overline{z}$-derivatives, and consequently
such kernels are not $\overline{z}-$smoothing of any positive order.
Therefore one needs to analyse the kernels of boundary order $\geq 0$ that
arise in the current setting more carefully in order to obtain the estimates
stated in the Main Theorem. It will be convenient to introduce the following
notation.

\begin{definition}
The symbol $\Gamma _{\lambda }$ denotes an admissible kernel of (boundary)
order $\geq \lambda $. We shall denote by $\Gamma _{0,1/2}^{\overline{z}}$
(resp. $\Gamma _{0,2/3}^{\overline{z}}$ ) an admissible kernel of order $%
\geq 0$ which is $\overline{z}$- smoothing of any order $\delta <1/2$ (resp. 
$\delta <2/3)$ and tangentially smoothing of any order $\delta <1/3$.
\end{definition}

Proposition 1 thus states that $\Omega _{q}(W^{\mathcal{L}})$ is a kernel of
type $\Gamma _{0,2/3}^{\overline{z}}$ according to this definition.
Similarly, we note that by Theorem 3 a kernel of type $\Gamma _{1}$ is
(better than) of type $\Gamma _{0,2/3}^{\overline{z}},$ and, in fact, is
smoothing of order $\delta <1/3$ in \emph{all} directions.

\section{The Integrals $\protect\int_{bD}\overline{\partial }f\wedge \Omega
_{q}(W^{\mathcal{L}},B)$ and $\protect\int_{bD}f\wedge \overline{\partial }%
_{z}\Omega _{q-1}(W^{\mathcal{L}},B)$}

We recall (see [10], [12], for example) that for $0\leq q\leq n-2$ the
transition kernels $\Omega _{q}(W^{\mathcal{L}},B)$ are defined by 
\begin{equation*}
\Omega _{q}(W^{\mathcal{L}},B)=(2\pi i)^{-n}\sum_{\mu
=0}^{n-q-2}\sum_{k=0}^{q}a_{\mu ,k,q}W^{\mathcal{L}}\wedge B\wedge (%
\overline{\partial }_{\zeta }W^{\mathcal{L}})^{\mu }\wedge
\end{equation*}%
\begin{equation}
\hspace{1in}\wedge (\overline{\partial }_{\zeta }B)^{n-q-2-\mu }\wedge (%
\overline{\partial }_{z}W^{\mathcal{L}})^{k}\wedge (\overline{\partial }%
_{z}B)^{q-k},  \label{transition}
\end{equation}%
where the coefficients $a_{\mu ,k,q}$ are certain rational numbers, while $%
\Omega _{n-1}(W^{\mathcal{L}},B)\equiv 0$. Again, it is enough to consider $%
\left\vert \zeta -z\right\vert \leq \varepsilon /2$, so that $W^{\mathcal{L}%
}=g/\Phi $. It then follows from (\ref{transition}) and standard results
about CF form that on $bD$ the form $\Omega _{q}(W^{\mathcal{L}},B)$ is
given by a sum of terms%
\begin{equation}
A_{q,\mu k}=\frac{a_{\mu ,k,q}}{(2\pi i)^{n}}\frac{g\wedge \partial \beta
\wedge (\overline{\partial }_{\zeta }g)^{\mu }\wedge (\overline{\partial }%
_{z}g)^{k}\wedge (\overline{\partial }_{\zeta }\partial \beta )^{n-q-2-\mu
}\wedge (\overline{\partial }_{z}\partial \beta )^{q-k}}{\Phi ^{1+\mu +k}~%
\mathcal{\beta }^{n-\mu -k-1}}\text{,}  \label{detail}
\end{equation}%
where $a_{\mu ,k,q}\in \mathbb{Q},$ $0\leq \mu \leq n-q-2$ and $0\leq k\leq
q $. As in case of the kernel $\Omega _{q}(W^{\mathcal{L}})$, it follows
that 
\begin{equation*}
A_{q,\mu k}(W^{\mathcal{L}},B)=\frac{\mathcal{L}[\mu ](\mathcal{E}%
_{1}^{\#})^{k}}{\Phi ^{1+\mu +k}}\frac{\mathcal{E}_{1}}{\mathcal{\beta }%
^{n-\mu -k-1}}.
\end{equation*}%
Consequently the kernels $A_{q,\mu k}(W^{\mathcal{L}},B)$ are admissible of
boundary order $\geq 1$. The integral $\int_{bD}\overline{\partial }f\wedge
\Omega _{q}(W^{\mathcal{L}},B)$ is therefore covered by part c) in Theorem
3; in particular, its kernel is smoothing of order $\delta <1/3.$

Next, one checks that 
\begin{eqnarray*}
\overline{\partial }_{z}A_{q,\mu k}(W^{\mathcal{L}},B) &=&\frac{\mathcal{L}%
[\mu ](\mathcal{E}_{1}^{\#})^{k}}{\Phi ^{1+\mu +k}}\frac{\mathcal{E}_{0}}{%
\mathcal{\beta }^{n-\mu -k-1}}+ \\
&&+\frac{\mathcal{L}[\mu ](\mathcal{E}_{1}^{\#})^{k}\mathcal{E}_{2}^{\#}}{%
\Phi ^{1+\mu +k+1}}\frac{\mathcal{E}_{1}}{\mathcal{\beta }^{n-\mu -k-1}}+%
\frac{\mathcal{L}[\mu ](\mathcal{E}_{1}^{\#})^{k}}{\Phi ^{1+\mu +k}}\frac{%
\mathcal{E}_{2}}{\mathcal{\beta }^{n-\mu -k}}.
\end{eqnarray*}%
This shows that the kernels $\overline{\partial }_{z}A_{q,\mu k}(W^{\mathcal{%
L}},B)$ are admissible, and one easily verifies that $\overline{\partial }%
_{z}A_{q,\mu k}$ is of boundary order $\geq 0.$

Recall that $0\leq \mu \leq n-q-2$ and $0\leq k\leq q$, so that $0\leq \mu
+k\leq n-2$. We shall first consider the case $\mu +k\geq 1$, which occurs
only when $n\geq 3$.

\begin{lemma}
Suppose that $\mu +k\geq 1$. Then $\overline{\partial }_{z}A_{q,\mu k}(W^{%
\mathcal{L}},B)$ is a kernel of type $\Gamma _{0,1/2}^{\overline{z}}$.
\end{lemma}

\emph{Proof. }We had already seen that $\overline{\partial }_{z}A_{q,\mu
k}(W^{\mathcal{L}},B)$ is of order $\geq 0$. Applying a derivative with
respect to $\overline{z_{j}}$ to a factor $1/\Phi $ in any of the summands
of $\overline{\partial }_{z}A_{q,\mu k}(W^{\mathcal{L}},B)$ results in a
term estimated by a kernel $\Gamma _{0}$ of order $\geq 0$ multiplied by $%
\left\vert \mathcal{E}_{2}^{\#}/\Phi \right\vert $, and since $\left\vert
\Phi \right\vert \gtrsim \left\vert r(z)\right\vert ^{1-\delta }\left\vert
\zeta -z\right\vert ^{3\delta }$, the factor $\left\vert \mathcal{E}%
_{2}^{\#}/\Phi \right\vert $ will be bounded by $\left\vert \zeta
-z\right\vert ^{\alpha }\left\vert r(z)\right\vert ^{-1+\delta }$ for some $%
\alpha >0$ as soon as $\delta <2/3$. By a) in Theorem 3, the kernel $%
\left\vert \zeta -z\right\vert ^{\alpha }\left\vert \Gamma _{0}\right\vert $
is integrable uniformly in $z$ as soon as $\alpha >0.$ For the other
differentiations, note that since $\mu +k\geq 1$, each summand in $\overline{%
\partial }_{z}A_{q,\mu k}$ has at least one factor $\mathcal{L}[1]/\Phi $ or 
$\mathcal{E}_{1}^{\#}/\Phi $ with weight $\geq -2$ in addition to the one
factor $1/\Phi $ which is counted with weight $\geq -1$ in the calculation
of the order of $\overline{\partial }_{z}A_{q,\mu k}$. Differentiating the
numerator of any such factor of weight $\geq -2$ results in a term $\mathcal{%
E}_{0}/\Phi $, which can be estimated by $\frac{1}{\left\vert \zeta
-z\right\vert ^{3\delta }}\frac{1}{\left\vert r(z)\right\vert ^{1-\delta }}$%
, where the factor $1/\left\vert \zeta -z\right\vert ^{3\delta }$ is of
weight $>-2$ for any $\delta <2/3$. If the differentiation is applied to any
of the remaining factors in each of the summands, the order of the kernel
decreases at most by $1$ without affecting any of the factors $\mathcal{L}%
[1]/\Phi $ or $\mathcal{E}_{1}^{\#}/\Phi $. In order to compensate for this
decrease, one must extract a factor $\mathcal{E}_{1}^{\#}$ from~such a
factor of weight $-2$, leaving a factor of weight $-2$ multiplied with a
suitable power $\left\vert r(z)\right\vert ^{-1+\delta }.$ This follows as
before for a factor $\mathcal{E}_{1}^{\#}/\Phi ,$ since $1/\Phi $ is
estimated by $1/(\left\vert \zeta -z\right\vert ^{2}\left\vert
r(z)\right\vert ^{1/3})$. For factors $\mathcal{L}[1]/\Phi $, note that
according to (\ref{tanest})---after introducing $z$-diagonalizing
coordinates---one may estimate $\left\vert \tan \overline{\partial }\partial
r/\Phi \right\vert $ by terms of the type 
\begin{equation*}
\frac{1}{\left\vert r(z)\right\vert +\left\vert \zeta _{l}-z_{l}\right\vert
^{2}}\lesssim \frac{1}{\left\vert \zeta _{l}-z_{l}\right\vert ^{2\delta
}\left\vert r(z)\right\vert ^{1-\delta }}\lesssim \lbrack \frac{1}{%
\left\vert \zeta _{l}-z_{l}\right\vert ^{2\delta +1}\left\vert
r(z)\right\vert ^{1-\delta }}]\mathcal{E}_{1}^{\#}
\end{equation*}%
for a suitable $l\leq n-1$. (See [14] for details.) Here the factor $%
1/\left\vert \zeta _{l}-z_{l}\right\vert ^{2\delta +1}$ is of weight $\geq
-2 $ for any $\delta <1/2$. Altogether we thus proved that $\overline{%
\partial }_{z}A_{q,\mu k}$ is $\overline{z}-$ smoothing of any order $\delta
<1/2$. Finally, if one considers a tangential derivative $L_{j}^{z}$, $j\leq
n-1$, the same arguments apply as long as $\delta <1/3$. $\blacksquare $

\begin{remark}
We note that this last argument restricts the order of $\overline{z}$%
-smoothing to $\delta <1/2$, while in all other previous instances one has
the stronger estimates of order $\delta <2/3$.
\end{remark}

\section{The Kernel $\overline{\partial _{z}}A_{q-1,00}$}

We are thus left with $\overline{\partial _{z}}A_{q-1,00}$. This is the
critical and most delicate case. Note that this kernel contains a term $%
\frac{\mathcal{E}_{0}}{\Phi \mathcal{\beta }^{n-1}}$ (of order $\geq 0$);
however, differentiation with respect to $\overline{z_{j}}$ will result,
among others, in a term $\mathcal{E}_{1}/(\Phi \mathcal{\beta }^{n})$ which
is estimated, at best, by $\left\vert \mathcal{E}_{1}/\mathcal{\beta }%
^{n}\right\vert \left\vert r(z)\right\vert ^{-1}$. We see that $\overline{%
\partial _{z}}A_{q-1,00}$ contains terms which are not $\overline{z}$%
-smoothing of any positive order $\delta >0$. In order to proceed we need to
identify these critical terms and exploit certain approximate symmetries in
analogy to the method introduced in [10].

We begin by applying Stokes' theorem to replace the integral 
\begin{equation*}
\int_{bD}f\wedge \overline{\partial _{z}}A_{q-1,00}
\end{equation*}%
by an integral over $D.$ For this purpose we first extend $A_{q-1,00}$%
---that is $W^{\mathcal{L}}$ and $B$---from the boundary into $D$ without
introducing any new singularities, as follows. By the estimate (\ref%
{estimate}) for $\left\vert \Phi \right\vert $, as long as $-\varepsilon
<r(\zeta )<0$, one has $\left\vert \Phi \right\vert \gtrsim \left\vert
r\zeta )\right\vert $. Choose a $C^{\infty }$ function $\varphi $ on $%
\overline{D}$ so that $\varphi (\zeta )\equiv 1$ for $-\varepsilon /2\leq
r(\zeta )$ and $\varphi (\zeta )\equiv 0$ for $r(\zeta )\leq -3\varepsilon
/4.$ We then define the $(0,1)$ form 
\begin{equation*}
\widehat{W^{\mathcal{L}}}(\zeta ,z)=\varphi (\zeta )W^{\mathcal{L}}(\zeta ,z)
\end{equation*}%
on $\overline{D}\times \overline{D}-\{(\zeta ,\zeta ):\zeta \in bD\}$, so
that $\widehat{W^{\mathcal{L}}}(\zeta ,z)=W^{\mathcal{L}}(\zeta ,z)$ for $%
\zeta \in bD$.

We also define 
\begin{equation*}
\widehat{B}(\zeta ,z)=\frac{\partial _{\zeta }\beta }{\mathcal{P}(\zeta ,z)}%
\text{, where }\mathcal{P}(\zeta ,z)=\beta (\zeta ,z)+\frac{2r\left( \zeta
\right) r(z)}{\left\Vert \partial r\left( \zeta \right) \right\Vert
\left\Vert \partial r\left( z\right) \right\Vert }
\end{equation*}%
on $\overline{D}\times \overline{D}-\{(\zeta ,\zeta ):\zeta \in bD\}.$ (Note
that $r\left( \zeta \right) r(z)>0$ for $\zeta ,z\in D$.) Clearly $\widehat{B%
}(\zeta ,z)=B(\zeta ,z)$ for $\zeta \in bD$. By replacing $W^{\mathcal{L}}$
with $\widehat{W^{\mathcal{L}}}$ and $B$ with $\widehat{B}$ in $A_{q-1,00}$
one can therefore assume that $A_{q-1,00}$ extends to $\overline{D}\times 
\overline{D}-\{(\zeta ,\zeta ):\zeta \in bD\}$ without any singularities.

It then follows that 
\begin{equation*}
\int_{bD}f\wedge \overline{\partial _{z}}A_{q-1,00}=\int_{D}\overline{%
\partial _{\zeta }}f\wedge \overline{\partial _{z}}A_{q-1,00}+(-1)^{q}%
\int_{D}f\wedge \overline{\partial _{\zeta }\,}\overline{\partial _{z}}%
A_{q-1,00}.
\end{equation*}

\begin{remark}
When one considers kernels that are integrated over $D$, the definition of
admissible kernels and of their weighted order needs to be modified
appropriately. First of all, the dimension of the domain of integration is
now $2n,$ which leads to an increase in order by one. Also, each factor $%
r(z) $ or $r(\zeta )$ in the numerator increases the order by at least one.
Furthermore, since both $r(\zeta )$ and $\func{Im}F^{(r)}(\zeta ,z)$ are
used as coordinates in a neighborhood of $z$ for $dist(z,bD)<\varepsilon ,$
the weighted order is adjusted to account for the fact that by estimate (\ref%
{estimate}), now up to \emph{two} factors $\Phi $ are counted with weight $%
\geq -1.$ (See Definition 4.2 in [10] for more details.) In particular, it
follows that the (extended) kernel $\overline{\partial _{z}}A_{q-1,00}$,
which is admissible of boundary order $\geq 0$, is admissible of order $\geq
1$ over $D$.
\end{remark}

It is straightforward to prove the analogous version of Theorem 3 for
kernels integrated over $D$. In particular, one then obtains the following
result.

\begin{lemma}
The operator $T_{q-1,00}:C_{(0,q+1)}(\overline{D})\rightarrow
C_{(0,q)}^{1}(D)$ defined by 
\begin{equation*}
T_{q-1,00}(\psi )=\int_{D}\psi \wedge \overline{\partial _{z}}A_{q-1,00}
\end{equation*}%
is smoothing of order $\delta $ for any $\delta <1/3$ and $\overline{z}-$%
smoothing of any order $\delta <1.$
\end{lemma}

Note that because of the term $K\left\vert \zeta -z\right\vert ^{3}$%
contained in $\Phi $, the kernel $\overline{\partial _{z}}A_{q-1,00}$ is
only of class $C^{1}$ jointly in $(\zeta ,z)$ near points where $\zeta =z$.

We are thus left with estimating integrals of the kernel $\overline{\partial
_{\zeta }\,}\overline{\partial _{z}}A_{q-1,00}=\overline{\partial _{z}}\,%
\overline{\partial _{\zeta }}A_{q-1,00}$ over $D$. This kernel is readily
seen to be of order $\geq 0$, but it contains terms that are not $\overline{z%
}-$smoothing of any order $\delta >0$.

Proceeding as in [10], we introduce the kernel $\overline{L_{q-1}}%
=(-1)^{q}\ast A_{q-1,00}$ for $1\leq q\leq n-1$, where $\ast $ is the Hodge
operator acting in the variable $\zeta $ with respect to the standard inner
product of forms in $\mathbb{C}^{n}$. Note that in [10] the definition of $%
L_{q-1}$ involved $\Omega _{q-1}(\widehat{W^{\mathcal{L}}},\widehat{B})$,
while here we only take those summands $A_{q,\mu k}$ with $\mu +k=0.$ Since $%
A_{q-1,00}=-\ast \overline{L_{q-1}}$, one then obtains 
\begin{equation*}
(-1)^{q}\overline{\partial _{\zeta }}A_{q-1,00}=\ast \ast \overline{\partial
_{\zeta }}A_{q-1,00}=\ast (\ast \overline{\partial _{\zeta }})(-\ast 
\overline{L_{q-1}})=\ast \overline{\vartheta _{\zeta }L_{q-1}},
\end{equation*}%
where $\vartheta _{\zeta }=-\ast \partial _{\zeta }\ast $ is the (formal)
adjoint of $\overline{\partial }$. It follows that 
\begin{equation*}
(-1)^{q}\int_{D}f\wedge \overline{\partial _{\zeta }}A_{q-1,00}=\int_{D}f%
\wedge \ast \overline{\vartheta _{\zeta }L_{q-1}}=(f,\vartheta _{\zeta
}L_{q-1})_{D}\text{,}
\end{equation*}%
where the inner product is taken by integrating the pointwise inner product
of forms over $D$. Since $\overline{\partial _{z}}$ commutes with $\ast
_{\zeta }$, one has%
\begin{equation*}
(-1)^{q}\int_{D}f\wedge \overline{\partial _{z}}\overline{\partial _{\zeta }}%
A_{q-1,00}=(f,\partial _{z}\vartheta _{\zeta }L_{q-1})_{D}\text{.}
\end{equation*}

Let us introduce the Hermitian transpose $K^{\ast }$ of a double form $%
K=K(\zeta ,z)$ by%
\begin{equation*}
K^{\ast }(\zeta ,z)=\overline{K(z,\zeta )}\text{.}
\end{equation*}%
Note that $K^{\ast }(\zeta ,z)$ is the kernel of the adjoint $T^{\ast }$ of
the operator $T:f\rightarrow (f,K(\cdot ,z))_{D}$, i.e., $T^{\ast
}(f)\rightarrow (f,K^{\ast }(\cdot ,z))_{D}$.

One now writes 
\begin{equation*}
(f,\partial _{z}\vartheta _{\zeta }L_{q-1})_{D}=(f,\partial _{z}\vartheta
_{\zeta }L_{q-1}-[\partial _{z}\vartheta _{\zeta }L_{q-1}]^{\ast
})_{D}+(f,[\partial _{z}\vartheta _{\zeta }L_{q-1}]^{\ast })_{D}\text{ .}
\end{equation*}%
Since $[\partial _{z}\vartheta _{\zeta }L_{q-1}]^{\ast }=\overline{\partial
_{\zeta }}\overline{\vartheta _{z}}L_{q-1}^{\ast }$, if $f\in dom(\overline{%
\partial }^{\ast })$ one may integrate by parts in the second inner product,
resulting in $(f,\overline{\partial _{\zeta }}\overline{\vartheta _{z}}%
L_{q-1}^{\ast })_{D}=$ $(\overline{\partial }^{\ast }f,\overline{\vartheta
_{z}}L_{q-1}^{\ast })_{D}$.

Expanding the definition of admissible kernels to allow for factors $\Phi
^{\ast }$, with corresponding definition of order, one verifies that the
kernel $\vartheta _{z}\overline{L_{q-1}^{\ast }}$is admissible of order $%
\geq 1$, and consequently it is smoothing of order $\delta <1/3,$ and $%
\overline{z}$-smoothing of order $\delta <1$.

\section{The Critical Singularities}

We shall now carefully examine $\partial _{z}\vartheta _{\zeta
}L_{q-1}-[\partial _{z}\vartheta _{\zeta }L_{q-1}]^{\ast }$ and verify that
there is a cancellation of critical terms, so that the conjugate of this
kernel is partially smoothing as required for the Main Theorem.

We shall use the standard orthonormal frame $\omega _{1},...,\omega _{n}$
for $(1,0)$ forms on a neighborhood $U$ of $P\in bD$, with $\omega
_{n}=\partial r(\zeta )/\left\Vert \partial r(\zeta )\right\Vert $. After
shrinking $\varepsilon ,$ we may assume that $B(P,\varepsilon )\subset U.$
As usual, we shall focus on estimating integrals for fixed $z\in D$ with $%
\left\vert z-P\right\vert <\varepsilon /2$, and integration over $\zeta \in
D\cap U$ with $r(\zeta )\geq -\varepsilon /2$ and $\left\vert \zeta
-z\right\vert <\varepsilon /2$. Let $L_{1},...,L_{n}$ be the corresponding
dual frame of $(1,0)$ vectors, acting in the $\zeta $ variables. If $%
L=\sum_{k=1}^{n}a_{k}(\zeta )\frac{\partial }{\partial \zeta _{k}}$, we
denote by $L^{z}=L_{\zeta }^{z}=\sum_{k=1}^{n}a_{k}(\zeta )\frac{\partial }{%
\partial z_{k}}$ the corresponding vector field acting in the $z$ variables.
A similar convention applies to $\omega _{j}^{z}$, which we shall denote by $%
\theta _{j}$

One then readily verifies the following equations.

\begin{enumerate}
\item $\overline{\partial }r(\zeta )=\overline{\omega _{n}}\left\Vert
\partial r(\zeta )\right\Vert ;$

\item $\partial \beta =\sum_{j=1}^{n}(L_{j}\beta )\omega _{j}$ and $%
\overline{\partial }\beta =\sum_{j=1}^{n}(\overline{L_{j}}\beta )\overline{%
\omega _{j}};$

\item $\overline{\partial }\partial \beta =2\sum_{j=1}^{n}\overline{\omega
_{j}}\wedge \omega _{j}+\mathcal{E}_{1};$

\item $\overline{\partial _{z}}\partial \beta =-2\sum_{j=1}^{n}\overline{%
\theta _{j}}\wedge \omega _{j}+\mathcal{E}_{1};$

\item $L_{j}^{z}\beta =-L_{j}\beta +\mathcal{E}_{2}$ and $L_{j}L_{k}\beta =%
\mathcal{E}_{1};$

\item $L_{j}\mathcal{P}=\mathcal{E}_{1}$ and $L_{j}^{z}\mathcal{P}=\mathcal{E%
}_{1}$ for $j<n.$
\end{enumerate}

Somewhat more delicate are the following two formulas. They are the analogs
of Lemma 5.9 and Lemma 5.35 in [10], with the differences due to the fact
that in the present paper the defining function is not normalized, as it is
restricted to a special form so that its level surfaces remain pseudoconvex.
The definition of the extension $\mathcal{P}$ has been modified accordingly.
Since both formulas require \emph{exact} identification of the leading
terms, we include the details of the proofs.

\begin{lemma}
$L_{n}^{z}\mathcal{P=-}\frac{2}{\left\Vert \partial r(\zeta )\right\Vert }%
\overline{\Phi }$ $+\mathcal{E}_{0}r(\zeta )r(z)+\mathcal{E}_{1}r(\zeta )+%
\mathcal{E}_{2}.$
\end{lemma}

\emph{Proof. }We fix $\zeta \in U.$ After a unitary change of coordinates in
the $\zeta $ variables, one may assume that $\partial r/\partial \zeta
_{j}(\zeta )=0$ for $j<n$ and $\partial r/\partial \zeta _{n}(\zeta )>0,$ so
that $\left\Vert \partial r(\zeta )\right\Vert =\partial r/\partial \zeta
_{n}(\zeta )\sqrt{2}$ and $(L_{n}^{z})_{z}=\sqrt{2}\frac{\partial }{\partial
z_{n}}+\mathcal{E}_{1}.$ Then $L_{n}^{z}r(z)=\sqrt{2}\frac{\partial r}{%
\partial z_{n}}(\zeta )+\mathcal{E}_{1}=\left\Vert \partial r(\zeta
)\right\Vert +\mathcal{E}_{1}$. In this coordinate system one has%
\begin{eqnarray*}
\sqrt{2}\overline{\Phi }(\zeta ,z) &=&\sqrt{2}\partial r/\partial \overline{%
\zeta _{n}}(\zeta )(\overline{\zeta _{n}-z_{n}})+\mathcal{E}_{2}-\sqrt{2}%
r(\zeta )= \\
&=&\left\Vert \partial r(\zeta )\right\Vert (\overline{\zeta _{n}-z_{n}})-%
\sqrt{2}r(\zeta )+\mathcal{E}_{2}\text{,}
\end{eqnarray*}%
and%
\begin{eqnarray*}
L_{n}^{z}\mathcal{P(\zeta },z) &=&-\sqrt{2}(\overline{\zeta _{n}-z_{n}})+%
\mathcal{E}_{2}+\frac{2r(\zeta )}{\left\Vert \partial r(\zeta )\right\Vert }%
\frac{L_{n}^{z}r(z)}{\left\Vert \partial r(z)\right\Vert }+\mathcal{E}%
_{0}r(\zeta )r(z)= \\
&=&-\sqrt{2}(\overline{\zeta _{n}-z_{n}})+\mathcal{E}_{2}+\frac{2r(\zeta )}{%
\left\Vert \partial r(\zeta )\right\Vert }[1+\mathcal{E}_{1}]+\mathcal{E}%
_{0}r(\zeta )r(z)= \\
&=&-\sqrt{2}(\overline{\zeta _{n}-z_{n}})+\frac{2r(\zeta )}{\left\Vert
\partial r(\zeta )\right\Vert }+\mathcal{E}_{2}+\mathcal{E}_{1}r(\zeta )+%
\mathcal{E}_{0}r(\zeta )r(z)= \\
&=&-\frac{\sqrt{2}}{\left\Vert \partial r(\zeta )\right\Vert }~\left[
\left\Vert \partial r(\zeta )\right\Vert (\overline{\zeta _{n}-z_{n}})-\sqrt{%
2}r(\zeta )\right] +\mathcal{E}_{2}+\mathcal{E}_{1}r(\zeta )+\mathcal{E}%
_{0}r(\zeta )r(z)\text{.}
\end{eqnarray*}%
The proof is completed by combining these two equations.$\blacksquare $

\begin{lemma}
$2P-\sum_{j=1}^{n-1}\left\vert L_{j}\beta \right\vert ^{2}=\frac{4}{%
\left\Vert \partial r(\zeta )\right\Vert \left\Vert \partial r(z)\right\Vert 
}\left\vert \Phi \right\vert ^{2}+\mathcal{E}_{3}+\mathcal{E}_{2}r(\zeta )$.
\end{lemma}

\emph{Proof.} Here we fix $z$, and after a unitary coordinate change in $%
\zeta $ we may assume that $L_{j}=\sqrt{2}\frac{\partial }{\partial \zeta
_{j}}+\mathcal{E}_{1}$, so that $L_{j}\beta =\sqrt{2}(\overline{\zeta
_{j}-z_{j}})+\mathcal{E}_{2}$. Hence $\sum_{j=1}^{n-1}\left\vert L_{j}\beta
\right\vert ^{2}=2\sum_{j=1}^{n-1}\left\vert \zeta _{j}-z_{j}\right\vert
^{2}+\mathcal{E}_{3}$, and therefore 
\begin{equation}
2P-\sum_{j=1}^{n-1}\left\vert L_{j}\beta \right\vert ^{2}=2\left\vert \zeta
_{n}-z_{n}\right\vert ^{2}+\frac{4r(\zeta )r(z)}{\left\Vert \partial r(\zeta
)\right\Vert \left\Vert \partial r(z)\right\Vert }+\mathcal{E}_{3}\text{.}
\label{twoP}
\end{equation}%
Furthermore,

\begin{eqnarray*}
\sqrt{2}\overline{\Phi }(\zeta ,z) &=&\sqrt{2}\partial r/\partial \overline{%
\zeta _{n}}(z)(\overline{\zeta _{n}-z_{n}})+\mathcal{E}_{2}-\sqrt{2}r(\zeta
)= \\
&=&\left\Vert \partial r(z)\right\Vert (\overline{\zeta _{n}-z_{n}})-\sqrt{2}%
r(\zeta )+\mathcal{E}_{2}\text{,}
\end{eqnarray*}%
It follows that 
\begin{eqnarray*}
2\left\vert \Phi (\zeta ,z)\right\vert ^{2} &=&\left\Vert \partial
r(z)\right\Vert (\overline{\zeta _{n}-z_{n}})\sqrt{2}\Phi -2r(\zeta )\Phi +%
\mathcal{E}_{2}\Phi = \\
&=&\left\Vert \partial r(z)\right\Vert ^{2}\left\vert \zeta
_{n}-z_{n}\right\vert ^{2}+\left\Vert \partial r(z)\right\Vert (\overline{%
\zeta _{n}-z_{n}})\left[ -\sqrt{2}r(\zeta )+\mathcal{E}_{2}\right] + \\
&&-2r(\zeta )\Phi +\mathcal{E}_{2}\Phi \\
&=&\left\Vert \partial r(z)\right\Vert ^{2}\left\vert \zeta
_{n}-z_{n}\right\vert ^{2}-\sqrt{2}r(\zeta )\left[ \left\Vert \partial
r(z)\right\Vert (\overline{\zeta _{n}-z_{n}})+\sqrt{2}\Phi \right] +\mathcal{%
E}_{3}+\mathcal{E}_{2}\Phi \text{.}
\end{eqnarray*}%
By equation (\ref{sym}) one has 
\begin{eqnarray*}
\sqrt{2}\Phi &=&\sqrt{2}\Phi ^{\ast }+\mathcal{E}_{3}=\left\Vert \partial
r(\zeta )\right\Vert (\overline{z_{n}-\zeta _{n}})-\sqrt{2}r(z)+\mathcal{E}%
_{2}= \\
&=&-\left\Vert \partial r(z)\right\Vert (\overline{\zeta _{n}-z_{n}})-\sqrt{2%
}r(z)+\mathcal{E}_{2}\text{,}
\end{eqnarray*}%
where we used that $\left\Vert \partial r(\zeta )\right\Vert =\left\Vert
\partial r(z)\right\Vert +\mathcal{E}_{1}$. Inserting this equation into the
previous one and using $\mathcal{E}_{2}\Phi =\mathcal{E}_{3}+\mathcal{E}%
_{2}r(\zeta ),$ results in 
\begin{eqnarray*}
2\left\vert \Phi (\zeta ,z)\right\vert ^{2} &=&\left\Vert \partial
r(z)\right\Vert ^{2}\left\vert \zeta _{n}-z_{n}\right\vert ^{2}-\sqrt{2}%
r(\zeta )(-\sqrt{2}r(z)+\mathcal{E}_{2})+\mathcal{E}_{3}+\mathcal{E}%
_{2}r(\zeta )= \\
&=&\left\Vert \partial r(\zeta )\right\Vert \left\Vert \partial
r(z)\right\Vert \left\vert \zeta _{n}-z_{n}\right\vert ^{2}+2r(\zeta )r(z)+%
\mathcal{E}_{3}+\mathcal{E}_{2}r(\zeta )= \\
&=&\frac{1}{2}\left\Vert \partial r(\zeta )\right\Vert \left\Vert \partial
r(z)\right\Vert \left[ 2\left\vert \zeta _{n}-z_{n}\right\vert ^{2}+\frac{%
4r(\zeta )r(z)}{\left\Vert \partial r(\zeta )\right\Vert \left\Vert \partial
r(z)\right\Vert }\right] +\mathcal{E}_{3}+\mathcal{E}_{2}r(\zeta )\text{.}
\end{eqnarray*}%
After inserting equation (\ref{twoP}) and rearranging, one obtains the
equation claimed in the Lemma.$\blacksquare $

We shall now identify precisely the kernel $L_{q-1}$ and the critical
highest order singularity of $\partial _{z}\vartheta _{\zeta }L_{q-1}$ with
respect to the standard frames introduced above. To simplify notation we
shall replace $q-1$ with $q$ and consider $A_{q,00}$ and $L_{q}$ for $0\leq
q\leq n-2$. The computations follow closely those for the case $\mu =0$ in
[10]; therefore we shall just state the relevant formulas, and provide more
details only where critical differences arise.

From (\ref{detail}) one sees that 
\begin{eqnarray*}
A_{q,00} &=&\frac{a_{q}}{(2\pi i)^{n}}\frac{g\wedge \partial \beta \wedge (%
\overline{\partial _{\zeta }}\partial \beta )^{n-q-2}\wedge (\overline{%
\partial _{z}}\partial \beta )^{q}}{\Phi \mathcal{P}^{n-1}}= \\
&=&\frac{a_{q}}{(2\pi i)^{n}}\frac{\partial r\wedge \partial \beta \wedge (%
\overline{\partial _{\zeta }}\partial \beta )^{n-q-2}\wedge (\overline{%
\partial _{z}}\partial \beta )^{q}+\mathcal{E}_{2}}{\Phi \mathcal{P}^{n-1}}.
\end{eqnarray*}%
Then 
\begin{eqnarray*}
\overline{L_{q}} &=&(-1)^{q+1}\ast A_{q,00}= \\
&=&\overline{C_{q}}+\frac{\mathcal{E}_{2}}{\Phi \mathcal{P}^{n-1}}\text{,}
\end{eqnarray*}%
where 
\begin{equation*}
\overline{C_{q}}=(-1)^{q+1}\frac{a_{q}}{(2\pi i)^{n}}\ast _{\zeta }\frac{%
\partial r\wedge \partial \beta \wedge (\overline{\partial _{\zeta }}%
\partial \beta )^{n-q-2}\wedge (\overline{\partial _{z}}\partial \beta )^{q}%
}{\Phi \mathcal{P}^{n-1}}.
\end{equation*}%
It follows that 
\begin{equation*}
\overline{\partial _{z}\vartheta _{\zeta }L_{q}}=\overline{\partial _{z}}%
\overline{\vartheta _{\zeta }}\overline{C_{q}}+\overline{\partial _{z}}%
\overline{\vartheta _{\zeta }}\frac{\mathcal{E}_{2}}{\Phi \mathcal{P}^{n-1}}%
\text{.}
\end{equation*}

\begin{remark}
Note that in contrast to [10], the kernels $L_{q}=C_{q}+\frac{\mathcal{E}_{2}%
}{\overline{\Phi }\mathcal{P}^{n-1}}$ and $\partial _{z}\vartheta _{\zeta
}L_{q}$ analysed here only involve the term corresponding to $\mu =0$ in
[10]. Since in this paper we shall be concerned with $\overline{z}-$%
smoothing, we need to consider the conjugates $\overline{L_{q}}$, $\overline{%
C_{q}},$ and $\overline{\partial _{z}\vartheta _{\zeta }L_{q}}$, which are
the kernels that appear in the integral $\int f\wedge \ast \overline{%
\partial _{z}\vartheta _{\zeta }L_{q}}=(f,\partial _{z}\vartheta _{\zeta
}L_{q})_{D}$.
\end{remark}

\begin{lemma}
$\overline{\partial _{z}}\overline{\vartheta _{\zeta }}\frac{\mathcal{E}_{2}%
}{\Phi \mathcal{P}^{n-1}}$ is of type $\Gamma _{1}$.
\end{lemma}

The proof of this lemma involves a straightforward verification. Note that $%
\frac{\mathcal{E}_{2}}{\Phi \mathcal{P}^{n-1}}$ is of order $\geq 3$.
Differentiation with respect to $\zeta $ reduces the order by $1$ only,
since after differentiating $1/\Phi $ the resulting factor $1/\Phi ^{2}$ has
weight $-2$. Similarly, since $\overline{\partial _{z}}\Phi =\mathcal{E}%
_{2}^{\#}$, subsequent application of $\overline{\partial _{z}}$ also
reduces the order by $1$ only.$\blacksquare $

Next we represent $\overline{C_{q}}$ in terms of the local orthonormal
frames. By utilizing the various formulas recalled above, it follows that 
\begin{equation*}
\overline{C_{q}}=\frac{\gamma _{q}\left\Vert \partial r(\zeta )\right\Vert }{%
i^{n}\Phi \mathcal{P}^{n-1}}\ast _{\zeta }\sum\limits_{\substack{ \left\vert
Q\right\vert =q  \\ j<n}}\omega _{n}\wedge (L_{j}\beta )\omega _{j}\wedge (%
\overline{\omega }\wedge \omega )^{J}\wedge \omega ^{Q}\wedge \overline{%
\theta }^{Q}+\frac{\mathcal{E}_{2}}{\Phi \mathcal{P}^{n-1}}\text{,}
\end{equation*}%
where $\gamma _{q}$ is a \emph{real} constant. The summation is over all
strictly increasing $q$ - tuples $Q$ with $n\notin Q$, over $j<n$ with $%
j\notin Q$, and $J$ is the ordered $n-q-2$ tuple complementary to $njQ$ in $%
\{1,...,n\}$. Since $\ast \lbrack \omega ^{njQ}\wedge (\overline{\omega }%
\wedge \omega )^{J}]=b_{q}i^{n}\omega ^{njQ}$, where $b_{q}$ is real, it
follows that 
\begin{equation*}
\overline{C_{q}}=\widetilde{\gamma _{q}}\left\Vert \partial r(\zeta
)\right\Vert \sum\limits_{\substack{ \left\vert Q\right\vert =q  \\ j<n}}%
\frac{L_{j}\beta }{\Phi \mathcal{P}^{n-1}}\omega ^{njQ}\wedge \overline{%
\theta }^{Q}+\frac{\mathcal{E}_{2}}{\Phi \mathcal{P}^{n-1}}\text{,}
\end{equation*}%
for another \emph{real} constant $\widetilde{\gamma _{q}}$. By using Lemma
11 it then follows that

\begin{equation*}
\overline{\partial _{z}}\overline{\vartheta _{\zeta }}\overline{C_{q}}=%
\widetilde{\gamma _{q}}\left\Vert \partial r(\zeta )\right\Vert \overline{%
\,\partial _{z}}\overline{\vartheta _{\zeta }}\left[ \sum\limits_{\left\vert
Q\right\vert =q,\,j<n}\frac{L_{j}\beta }{\Phi \mathcal{P}^{n-1}}\omega
^{njQ}\wedge \overline{\theta }^{Q}\right] +\Gamma _{1}\text{.}
\end{equation*}

Let us introduce 
\begin{equation*}
\overline{C_{q}^{0}}=\sum\limits_{\left\vert Q\right\vert =q,\,j<n}\frac{%
L_{j}\beta }{\Phi \mathcal{P}^{n-1}}\omega ^{njQ}\wedge \overline{\theta }%
^{Q}\text{.}
\end{equation*}%
The following theorem contains the heart of the analysis of $\partial
_{z}\vartheta _{\zeta }L_{q}-[\partial _{z}\vartheta _{\zeta }L_{q}]^{\ast
}. $

\begin{theorem}
For $0\leq q\leq n-2$ the kernel 
\begin{equation*}
\overline{\partial _{z}}\overline{\vartheta _{\zeta }}\overline{C_{q}^{0}}-[%
\overline{\partial _{z}}\overline{\vartheta _{\zeta }}\overline{C_{q}^{0}}%
]^{\ast }
\end{equation*}%
is of type $\Gamma _{0,2/3}^{\overline{z}}$.
\end{theorem}

\begin{corollary}
The operator 
\begin{equation*}
f\rightarrow (f,\partial _{z}\vartheta _{\zeta }L_{q-1}-[\partial
_{z}\vartheta _{\zeta }L_{q-1}]^{\ast })_{D}
\end{equation*}%
is $\overline{z}-$smoothing of order $\delta <2/3$ and tangentially
smoothing of order $\delta <1/3$ for $1\leq q\leq n-1.$
\end{corollary}

This follows from the Theorem by using Lemma 12 and also by observing that
differentiation of $\left\Vert \partial r(\zeta )\right\Vert $ results in an
error term of type $\Gamma _{1}$. Similarly, when considering the difference 
$[...]-[...]^{\ast }$, the substitution $\left\Vert \partial r(\zeta
)\right\Vert =\left\Vert \partial r(z)\right\Vert +$ $\mathcal{E}_{1}$ leads
to an error term of the same type.$\blacksquare $

\begin{remark}
In order to be consistent with notations and formulas in [10], in the proof
of the theorem we shall analyse $\Delta _{q}=$ $\partial _{z}\vartheta
_{\zeta }C_{q}^{0}-[\partial _{z}\vartheta _{\zeta }C_{q}^{0}]^{\ast }$. In
the end we must verify that its \textbf{conjugate} $\overline{\Delta _{q}}$
is of type $\Gamma _{0,2/3}^{\overline{z}}$.
\end{remark}

Since $C_{q}^{0}$ is a double form of type $(0,q+2)$ in $\zeta $ and type $%
(q,0)$ in $z$, the form $\partial _{z}\vartheta _{\zeta }C_{q}^{0}$ is of
type $(0,q+1)$ in $\zeta $ and type $(q+1,0)$ in $z$. Consequently, 
\begin{equation*}
\partial _{z}\vartheta _{\zeta }C_{q}^{0}=\vartheta _{\zeta }\partial
_{z}C_{q}^{0}=\sum_{\left\vert L\right\vert =q+1}\left( \sum_{\left\vert
K\right\vert =q+1}A_{KL}\overline{\omega }^{K}\right) \wedge \theta ^{L},
\end{equation*}%
where the sums are taken over all strictly increasing $q+1$ tuples $L$ and $%
K.$ It follows that 
\begin{equation*}
\partial _{z}\vartheta _{\zeta }C_{q}^{0}-[\partial _{z}\vartheta _{\zeta
}C_{q}^{0}]^{\ast }=\sum_{\left\vert L\right\vert =q+1}\left(
\sum_{\left\vert K\right\vert =q+1}[A_{KL}-(A_{LK})^{\ast }]\overline{\omega 
}^{K}\right) \wedge \theta ^{L}\text{.}
\end{equation*}

In the next section we shall identify the coefficients $A_{KL}$ precisely in
order to verify that $[A_{KL}-(A_{LK})^{\ast }]$ is of order $\geq 0$ and
that its \emph{conjugate} is at least of type $\Gamma _{0,2/3}^{\overline{z}%
} $.

\section{The Approximate Symmetries}

The computation of $\vartheta _{\zeta }\partial _{z}C_{q}^{0}$ uses the
expressions for $\partial _{z}$ and $\vartheta _{\zeta }$ in terms of the
standard adapted boundary frames plus error terms which do not involve
differentiation. These error terms result---in the end---in kernels which
are conjugates of admissible kernels of order $\geq 1$, and hence will be
ignored in the discussion that follows.

As usual $\varepsilon _{lQ}^{L}$ denotes the sign of the permutation which
carries the ordered $(q+1)$ -tuple $lQ$ into the ordered $(q+1)$ -tuple $L$
if $lQ=L$ as sets, and $\varepsilon _{lQ}^{L}=0$ otherwise. Let us introduce 
$m_{lj}=\frac{1}{\overline{\Phi }}L_{l}^{z}(\frac{\overline{L_{j}}\beta }{%
\mathcal{P}^{n-1}})$ for $1\leq j<n$ and $1\leq l\leq n$.

\begin{lemma}
For any $K,L$ one has 
\begin{equation*}
A_{KL}=-\sum\limits_{\substack{ Q  \\ j,l,k}}\varepsilon
_{lQ}^{L}\varepsilon _{kK}^{njQ}(L_{k}m_{lj})+\overline{\Gamma _{0,2/3}^{%
\overline{z}}}+\Gamma _{1}\text{.}
\end{equation*}
\end{lemma}

Proof. This Lemma is the analogue of Lemma 5.6 in [10]. In the present case,
the computation shows that the error term is of the form $\mathcal{E}_{3}/[%
\overline{\Phi }^{3}\mathcal{P}^{n-1}]+\Gamma _{1}$, where the first term is
only of order $0$. (In case $D$ is strictly pseudoconvex, and hence $%
\left\vert \Phi \right\vert \gtrsim \left\vert \zeta -z\right\vert ^{2}$,
this term is $\Gamma _{1}$ as well.) However, one readily checks that its
conjugate $\mathcal{E}_{3}/[\Phi ^{3}\mathcal{P}^{n-1}]$ is of type $\Gamma
_{0,2/3}^{\overline{z}}$. $\blacksquare $

We are therefore left with 
\begin{equation}
A_{KL}^{(0)}=-\sum\limits_{\substack{ Q  \\ j,l,k}}\varepsilon
_{lQ}^{L}\varepsilon _{kK}^{njQ}(L_{k}m_{lj}).  \label{AKL}
\end{equation}

Only terms with $j<n$ appear with nonzero coefficients. In the following it
will be assumed that $j<n$.

For $l<n$ one has%
\begin{equation}
m_{lj}=\frac{-2\delta _{lj}}{\overline{\Phi }\mathcal{P}^{n-1}}-(n-1)\frac{%
L_{l}^{z}\beta \overline{L_{j}}\beta }{\overline{\Phi }\mathcal{P}^{n}}+%
\frac{\mathcal{E}_{1}}{\overline{\Phi }\mathcal{P}^{n-1}}\text{,}
\label{Mlj}
\end{equation}%
while by using Lemma 9 one obtains%
\begin{eqnarray}
m_{nj} &=&\frac{1}{\overline{\Phi }}[\frac{L_{n}^{z}\overline{L_{j}}\beta }{%
\mathcal{P}^{n-1}}-(n-1)\frac{(L_{n}^{z}\mathcal{P)}\overline{L_{j}}\beta }{%
\mathcal{P}^{n}}]=  \label{Mnj} \\
&=&\frac{2(n-1)}{\left\Vert \partial r(\zeta )\right\Vert }\frac{\overline{%
L_{j}}\beta }{\mathcal{P}^{n}}+\frac{\mathcal{E}_{1}}{\overline{\Phi }%
\mathcal{P}^{n-1}}+\frac{\mathcal{E}_{1}r(\zeta )r(z)+\mathcal{E}_{2}r(\zeta
)+\mathcal{E}_{3}}{\overline{\Phi }\mathcal{P}^{n}}.  \notag
\end{eqnarray}%
Note that all the error terms are of order $\geq 2$, and that they have only 
\emph{one} factor $\Phi $ in the denominator. Consequently, applying $L_{k}$
to the error terms results in $\Gamma _{1}$ terms. So only the leading terms
of $m_{lj}$ identified above need to be considered in the following analysis.

We now come to the proof of Theorem 13. Since all relevant error terms are
of type $\Gamma _{0,2/3}^{\overline{z}}$ or better, it is enough to examine $%
\overline{A_{KL}^{(0)}}-\overline{A_{LK}^{(0)\ast }}$, where 
\begin{equation*}
\overline{A_{KL}^{(0)}}=-\sum\limits_{\substack{ Q  \\ j,l,k}}\varepsilon
_{L}^{lQ}\varepsilon _{kK}^{njQ}(\overline{L_{k}m_{lj}}).
\end{equation*}

We need to consider separate cases, depending on whether $n$ is in $K$,
resp. $L$, or not.

\textbf{Case 1. }$n\in K$ and $n\in L$.

In this case the computations in [10] apply without further changes, subject
to the adjustments due to the fact that the defining function $r$ is not
normalized. Combined with $\left\Vert \partial r(\zeta )\right\Vert
=\left\Vert \partial r(z)\right\Vert +\mathcal{E}_{1},$ it follows that 
\begin{equation*}
\overline{A_{KL}^{(0)}-A_{LK}^{(0)\ast }}=\Gamma _{1}.
\end{equation*}

\textbf{Case 2.} $n\notin K$ and $n\notin L$.

In this case $\varepsilon _{kK}^{njQ}\neq 0$ only for $k=n.$ Hence%
\begin{equation*}
A_{KL}^{(0)}=-\sum\limits_{\substack{ Q  \\ j,l}}\varepsilon
_{lQ}^{L}\varepsilon _{nK}^{njQ}(L_{n}m_{lj})+\Gamma _{1}=-\sum\limits 
_{\substack{ l\in L  \\ j\in K}}\varepsilon _{lQ}^{L}\varepsilon
_{K}^{jQ}(L_{n}m_{lj})+\Gamma _{1}.
\end{equation*}%
Lemma 5.21 in [10] needs to be replaced by

\begin{lemma}
$L_{n}m_{lj}-(L_{n}m_{jl})^{\ast }=\overline{\Gamma _{0,2/3}^{\overline{z}}}$
for $j,l<n.$
\end{lemma}

Assuming the lemma, one obtains---after replacing $j$ with $l$ in the last
equation--- 
\begin{eqnarray*}
A_{LK}^{(0)\ast } &=&-\sum\limits_{j,l<n}\varepsilon _{lQ}^{K}\varepsilon
_{L}^{jQ}(L_{n}m_{lj})^{\ast }+\Gamma _{1}=-\sum\limits_{\substack{ j\in L 
\\ l\in K}}\varepsilon _{lQ}^{K}\varepsilon _{L}^{jQ}(L_{n}m_{jl})+\overline{%
\Gamma _{0,2/3}^{\overline{z}}}= \\
&=&A_{KL}^{(0)}+\overline{\Gamma _{2/3}^{\overline{z}}}\text{.}
\end{eqnarray*}

\bigskip

To prove the Lemma, one uses equation (\ref{Mlj}). The calculation of $%
L_{n}m_{lj}$ proceeds as in [10] with the obvious changes. By using $%
\left\Vert \partial r(\zeta )\right\Vert =\left\Vert \partial
r(z)\right\Vert +\mathcal{E}_{1}$ one obtains 
\begin{eqnarray*}
(L_{n}m_{jl})^{\ast }-L_{n}m_{lj} &=&-\frac{2\delta _{lj}\left\Vert \partial
r(\zeta )\right\Vert }{\mathcal{P}^{n-1}}\left[ \frac{1}{\overline{\Phi
^{\ast 2}}}-\frac{1}{\overline{\Phi ^{2}}}\right] -4\frac{\delta _{lj}(n-1)}{%
\left\Vert \partial r(\zeta )\right\Vert \mathcal{P}^{n}}\left[ \frac{\Phi
^{\ast }}{\overline{\Phi ^{\ast }}}-\frac{\Phi }{\overline{\Phi }}\right] +
\\
&&-(n-1)\frac{L_{j}^{z}\beta \overline{L_{l}}\beta \left\Vert \partial
r(\zeta )\right\Vert }{\mathcal{P}^{n}}\left[ \frac{1}{\overline{\Phi ^{\ast
2}}}-\frac{1}{\overline{\Phi ^{2}}}\right] + \\
&&-\frac{2n(n-1)}{\left\Vert \partial r(\zeta )\right\Vert }\frac{%
L_{j}^{z}\beta \overline{L_{l}\beta }}{\mathcal{P}^{n+1}}\left[ \frac{\Phi
^{\ast }}{\overline{\Phi ^{\ast }}}-\frac{\Phi }{\overline{\Phi }}\right]
+\Gamma _{1}.
\end{eqnarray*}%
In the strictly pseudoconvex case the differences in $[...]$ are of higher
order than the terms individually, resulting in $(L_{n}m_{jl})^{\ast
}-L_{n}m_{lj}=\Gamma _{1}.$ In the present case, only a weaker result holds,
as follows. Note that---after taking conjugates---one has 
\begin{equation*}
\frac{1}{\Phi ^{\ast 2}}-\frac{1}{\Phi ^{2}}=\frac{(\Phi -\Phi ^{\ast
})(\Phi +\Phi ^{\ast })}{\Phi ^{\ast 2}\Phi ^{2}}=\frac{\mathcal{E}_{3}}{%
\Phi ^{\ast 2}\Phi }+\frac{\mathcal{E}_{3}}{\Phi ^{\ast }\Phi ^{2}}\text{,}
\end{equation*}%
where we have used the approximate symmetry (\ref{sym}) in the second
equation. It now readily follows that $\mathcal{E}_{3}/[\mathcal{P}%
^{n-1}\Phi ^{\ast 2}\Phi ]$ and $\mathcal{E}_{3}/[\mathcal{P}^{n-1}\Phi
^{\ast }\Phi ^{2}]$ (while of order $0,$ and hence not smoothing as in the
strictly pseudoconvex case) are in fact of type $\Gamma _{0,2/3}^{\overline{z%
}}$. Since $\mathcal{E}_{2}/\mathcal{P}^{n}$ is estimated by $\mathcal{E}%
_{0}/\mathcal{P}^{n-1}$, the same argument works for the third term above.
For the conjugate of the second term note that 
\begin{eqnarray*}
\frac{1}{\mathcal{P}^{n}}\left[ \frac{\overline{\Phi ^{\ast }}}{\Phi ^{\ast }%
}-\frac{\overline{\Phi }}{\Phi }\right] &=&\frac{1}{\mathcal{P}^{n}}\left[ 
\frac{\Phi \overline{\Phi ^{\ast }}-\Phi ^{\ast }\overline{\Phi }}{\Phi
^{\ast }\Phi }\right] = \\
&=&\frac{1}{\mathcal{P}^{n}}\frac{\mathcal{E}_{3}}{\Phi ^{\ast }\Phi }%
=\Gamma _{1}\text{.}
\end{eqnarray*}%
The fourth term is estimated the same way by first estimating $\mathcal{E}%
_{2}/\mathcal{P}^{n+1}$ by $\mathcal{E}_{0}/\mathcal{P}^{n}$.

\textbf{Case 3. }The mixed case $n\in K$ and $n\notin L$.

As in [10], this is---computationally---the most complicated case. On the
other hand, aside from the differences as noted, for example, in Lemmas 9
and 10, the details of the proof essentially carry over from [10] to the
case considered here, with the result that one has $\overline{A_{KL}^{(0)}}-%
\overline{A_{LK}^{(0)\ast }}=\Gamma _{1}.$ In some more detail, since $n\in
K $ there is exactly one ordered $q$-tuple $J$ such that $K=J\cup \{n\}$,
and one then has $\varepsilon _{nJ}^{K}A_{KL}^{(0)}=A_{(nJ)L}^{(0)}$. Note
that we need to identify the leading terms of both $\overline{A_{KL}^{(0)}}$
and $\overline{A_{LK}^{(0)}}$. Let us first consider the simpler term $%
\overline{A_{LK}^{(0)}}$. After interchanging $L$ and $K$ in equation (\ref%
{AKL}), one has%
\begin{equation*}
A_{LK}^{(0)}=-\sum\limits_{\substack{ Q  \\ j,l,k}}\varepsilon
_{lQ}^{K}\varepsilon _{kL}^{njQ}(L_{k}m_{lj})+\Gamma _{1}\text{.}
\end{equation*}%
Since $n\notin L,$ the factor $\varepsilon _{lQ}^{K}\varepsilon
_{kL}^{njQ}\neq 0$ only if $k=n$ and $l=n,$ and furthermore $Q=J$. Therefore
the leading term of $A_{LK}^{(0)}$, i.e. the sum, is different from $0$ only
if $J\subset L$ so that $\varepsilon _{L}^{jJ}\neq 0$ only for that unique $%
j $ for which $L=J\cup \{j\}.$ It follows that for $j<n$ one has 
\begin{equation*}
A_{(jJ)(nJ)}^{(0)}=\varepsilon _{L}^{jJ}\varepsilon
_{nJ}^{K}A_{LK}^{(0)}+\Gamma _{1}=-L_{n}m_{nj}+\Gamma _{1}\text{.}
\end{equation*}%
Since $L_{n}\mathcal{P=(}\overline{L_{n}^{z}\mathcal{P}}\mathcal{)}^{\ast }$%
, Lemma 9 implies that $L_{n}\mathcal{P=-}\frac{2}{\left\Vert \partial
r(z)\right\Vert }\Phi ^{\ast }$ $+\mathcal{E}_{0}r(\zeta )r(z)+\mathcal{E}%
_{1}r(z)+\mathcal{E}_{2}$. By using this equation and (\ref{Mnj}) it follows
that%
\begin{equation}
A_{(jJ)(nJ)}^{(0)}=-\frac{4n(n-1)}{\left\Vert \partial r(\zeta )\right\Vert
\left\Vert \partial r(z)\right\Vert }\frac{\overline{L_{j}\beta }\,\Phi
^{\ast }}{\mathcal{P}^{n+1}}+\Gamma _{1}.  \label{A jJ nJ}
\end{equation}

Finally we calculate $\overline{A_{KL}^{(0)}}.$ With $J$ as before, one has 
\begin{equation*}
A_{KL}^{(0)}=-\sum\limits_{\substack{ Q  \\ j,l,k}}\varepsilon
_{lQ}^{L}\varepsilon _{kK}^{njQ}(L_{k}m_{lj})+\Gamma _{1}=\varepsilon
_{K}^{nJ}\sum\limits_{\substack{ n\notin Q  \\ j,l,k<n}}\varepsilon
_{lQ}^{L}\varepsilon _{kJ}^{jQ}(L_{k}m_{lj})+\Gamma _{1}\text{.}
\end{equation*}%
Continuing with the intricate calculations as in case Id) on pp. 237-239 in
[10] and using Lemma 10 above in place of Lemma 5.35 in [10], one obtains%
\begin{equation*}
A_{KL}^{(0)}=-\varepsilon _{K}^{nJ}\frac{4n(n-1)}{\left\Vert \partial
r(\zeta )\right\Vert \left\Vert \partial r(z)\right\Vert }%
\sum_{l<n}\varepsilon _{lJ}^{L}\frac{(L_{l}^{z}\beta )\Phi }{\mathcal{P}%
^{n+1}}+\Gamma _{1}\text{.}
\end{equation*}%
Here the only nonzero term in the sum arises for that unique $l,$ for which $%
L=J\cup \{l\}.$ Consequently, the last formula implies

\begin{equation}
A_{(nJ)(lJ)}^{(0)}=-\frac{4n(n-1)}{\left\Vert \partial r(\zeta )\right\Vert
\left\Vert \partial r(z)\right\Vert }\frac{(L_{l}^{z}\beta )\Phi }{\mathcal{P%
}^{n+1}}+\Gamma _{1}\text{.}  \label{A nJ lJ}
\end{equation}

Let us now consider $\overline{A_{KL}^{(0)}}-\overline{A_{LK}^{(0)\ast }}$.
As the preceding formulas show, each summand is of type $\Gamma _{1}$ except
in the case that for the unique $q$-tuple $J\subset \{1,2,...,n-1\}$ with $%
K=J\cup \{n\}$, $L$ satisfies $L=J\cup \{l\}$ as sets for some unique $l<n$.
In this latter case, equations (\ref{A nJ lJ}) and (\ref{A jJ nJ}) imply 
\begin{equation*}
\overline{A_{(nJ)(lJ)}^{(0)}}-\overline{A_{(lJ)(nJ)}^{(0)\ast }}=-\frac{%
4n(n-1)}{\left\Vert \partial r(\zeta )\right\Vert \left\Vert \partial
r(z)\right\Vert }\frac{1}{\mathcal{P}^{n+1}}\overline{\left[ (L_{l}^{z}\beta
)\Phi \mathcal{-(}\overline{L_{l}\beta }\,\mathcal{)}^{\ast }\Phi \right] }%
+\Gamma _{1}=\Gamma _{1}\text{.}
\end{equation*}%
The last equation holds because $\,\mathcal{(}\overline{L_{l}\beta }\,%
\mathcal{)}^{\ast }=L_{l}^{z}\beta $.

\textbf{Case 4.} $n\notin K$ and $n\in L.$

This is reduced to Case 3 by noting that $\overline{A_{KL}^{(0)}}-\overline{%
A_{LK}^{(0)\ast }}=-$ $(\overline{A_{LK}^{(0)}}-\overline{A_{KL}^{(0)\ast }}%
)^{\ast }.$

It thus follows that for all $K$ and $L$ one has $\overline{A_{KL}^{(0)}}-%
\overline{A_{LK}^{(0)\ast }}=\Gamma _{0,2/3}^{\overline{z}}$. This completes
the proof of Theorem 13.$\blacksquare $

\section{Proof of the Main Theorem.}

In sections 4 - 8 we have analysed the integrals that appear in the
representation (\ref{Kop}) of the boundary operator $S^{bD}$. By combining
these results it follows that for all $f\in \mathcal{D}_{qU}^{1}$ and $z\in
D\cap U$ one has the estimates 
\begin{equation*}
\left\vert \overline{L_{j}^{z}}S^{bD}(f)(z)\right\vert \leq C_{\delta }\cdot
dist(z,bD)^{\delta -1}\cdot Q_{0}(f)\emph{\ }\text{\emph{for}}\emph{\ }%
j=1,...,n\emph{\ }\text{\emph{and any}}\emph{\ }\delta <1/2
\end{equation*}%
and 
\begin{equation*}
\left\vert L_{j}^{z}S^{bD}(f)(z)\right\vert \leq C_{\delta }\cdot
dist(z,bD)^{\delta -1}\cdot Q_{0}(f)\emph{\ }\text{\emph{for}}\emph{\ }%
j=1,...,n-1\emph{\ }\text{\emph{and any}}\emph{\ }\delta <1/3.
\end{equation*}

We have thus completed the proof of parts 2) and 3) of the Main Theorem. As
noted earlier, part 1) follows trivially from the classical estimate (\ref%
{iso}) for $S^{iso}$.

Finally we prove the statement about the normal components of $f$. We shall
use the following Lemma, which is a routine variation of classical estimates
for the BMK kernel. For $0\leq \alpha <1$ define $C_{(0,1)}^{-\alpha
}(D)=\{g\in C_{(0,1)}(D):$ $\sup_{z\in D}\left\vert g(z)\right\vert
dist(z,bD)^{\alpha }<\infty \}$, with the norm $\left\vert g\right\vert
_{-\alpha }$ defined by the relevant supremum.

\begin{lemma}
The operator $T^{BM}:C_{(0,1)}^{-\alpha }(D)\rightarrow C(D)$ defined by 
\begin{equation*}
T^{BM}(g)=\int_{D}g(\zeta )\wedge \Omega _{0}(B)
\end{equation*}%
satisfies the estimate%
\begin{equation}
\left\vert T^{BM}(g)\right\vert _{1-\alpha ^{\prime }}\lesssim \left\vert
g\right\vert _{-\alpha }\text{ for any }\alpha ^{\prime }>\alpha .
\label{wH}
\end{equation}
\end{lemma}

Now suppose $f\in \mathcal{D}_{qU}^{1}$ and let $f_{J}$ be a normal
component of $f$, so that $\left. f_{J}\right\vert _{bD}=0$. Decompose $%
f_{J}=h_{J}+[S^{iso}(f)]_{J}$, where $h=S^{bD}(f).$ We already know by
estimate (\ref{iso}) that $\left\vert S^{iso}(f)\right\vert _{\alpha
}\lesssim Q_{0}(f)$ for any $\alpha <1$. Note that on $bD\cap U$ one has $%
h_{J}=-[S^{iso}(f)]_{J}$, so that $\left\vert (\left. h_{J}\right\vert
_{bD\cap U})\right\vert _{\alpha }\lesssim \left\vert S^{iso}(f)\right\vert
_{\alpha }\lesssim Q_{0}(f)$ as well. By standard properties of the BM
kernel $\Omega _{0}(B)$ it follows that $\int_{bD}h_{J}\,\Omega _{0}(B)$
satisfies the same estimate on $\overline{D}\cap U$ if $\alpha >0$. By the
case $q=0$ of the BMK representation formula (\ref{BMK}) applied to $h_{J}$
one has 
\begin{equation}
h_{J}=\int_{bD}h_{J}\,\Omega _{0}(B)-\int_{D}\overline{\partial }h_{J}\wedge
\Omega _{0}(B).  \label{BMn}
\end{equation}%
Given $\delta <1/2$, choose $\delta ^{\prime }$ with $\delta <\delta
^{\prime }<1/2.$ By part 2) of the Main Theorem, $\overline{\partial }%
h_{J}\in C_{(0,1)}^{-(1-\delta ^{\prime })}(D)$, with $\left\vert \overline{%
\partial }h_{J}\right\vert _{-(1-\delta ^{\prime })}\leq C_{\delta ^{\prime
}}Q_{0}(f)$. It then follows from Lemma 18 that 
\begin{equation*}
\left\vert T^{BM}(\overline{\partial }h_{J})\right\vert _{\delta }\lesssim
Q_{0}(f).
\end{equation*}

Each summand in the representation (\ref{BMn}) therefore satisfies the
desired H\"{o}lder estimate, so that 
\begin{equation*}
\left\vert h_{J}\right\vert _{\Lambda ^{\delta }(\overline{D}\cap
U)}\lesssim Q_{0}(f)\text{. }
\end{equation*}%
Since $f_{J}=h_{J}+[S^{iso}(f)]_{J}$, the required estimate $\left\vert
f_{J}\right\vert _{\delta }\lesssim Q_{0}(f)$ holds as well. $\blacksquare 
\vspace{0.5in}$\newpage

\textbf{References}

[1]\qquad Catlin, D.: Subelliptic estimates for the $\overline{\partial }-$%
Neumann problem on pseudoconvex domains. \textit{Ann. Math. }\textbf{126 }%
(1987), 131- 191.

[2]\qquad Cumenge, A.: Estim\'{e}es Lipschitz optimales dans les convexes de
type fini. \textit{C. R. Acad. Sci. Paris} \textbf{325} (1997), 1077-1080.

[3]\qquad D'Angelo, J. F.: Real hypersurfaces, orders of contact, and
applications. \textit{Ann. Math.} \textbf{115} (1982), 615-637.

[4]\qquad Diederich, K., and Fornaess, J. E.: Pseudoconvex domains with
real-analytic boundary. \textit{Ann. Math.} \textbf{107} (1978), 371-384.

[5]\qquad Diederich, K., Fischer, B., and Fornaess, J. E.: H\"{o}lder
estimates on convex domains of finite type. \textit{Math. Z.} \textbf{232}
(1999), 43 - 61.

[6]\qquad Folland, G., and Kohn, J.J.: \textit{The Neumann Problem for the
Cauchy-Riemann Complex.} Princeton Univ. Press, 1972.

[7]\qquad Kohn, J.J.: Boundary behavior of $\overline{\partial }$ on weakly
pseudoconvex manifolds of dimension two. \textit{J. Diff. Geometry} \textbf{6%
} (1972), 523 -542.

[8]\qquad Kohn, J.J.: Subelllipticity of the $\overline{\partial }-$Neumann
Problem on Pseudoconvex Domains: Sufficient Conditions. \textit{Acta math.} 
\textbf{142} (1979), 79-122.

[9]\qquad Kohn, J.J., and Nirenberg, L.: A pseudoconvex domain not admitting
a holomorphic support function. \textit{Math. Ann. }\textbf{201 }(1973), 265
- 268.

[10]\qquad Lieb, I., and Range, R. M.: On integral representations and
a-priori Lipschitz estimates for the canonical solution of the $\overline{%
\partial }$-equation. Math. Ann. \textbf{265} (1983), 221 - 251.

[11]\qquad Lieb, I., and Range, R. M.: Integral representations and
estimates in the theory of the $\overline{\partial }-$Neumann problem. Ann.
Math. \textbf{123} (1986), 265 - 301.

[12]\qquad Range, R. M.: \textit{Holomorphic Functions and Integral
Representations in Several Complex Variables.} Springer Verlag New York
1986, corrected 2nd. printing 1998.

[13]\qquad Range, R. M.: A pointwise basic estimate and H\"{o}lder
multipliers for the $\overline{\partial }-$Neumann problem on pseudoconvex
domains. \textit{Preliminary Report.} \textit{arXiv:1106.3132 (June 2011)}

[14]\qquad Range, R. M.: An integral kernel for weakly pseudoconvex domains.
Math. Ann. (2012), DOI 10.1007/s00208-012-0863-4.

[15]\qquad Siu, Y.-T.: Effective termination of Kohn's algorithm for
subelliptic multipliers. \textit{Pure Appl. Math. Quarterly} \textbf{6}
(2010), no. 4, Special Issue: In honor of Joseph J. Kohn.\vspace{0.5in}

R. Michael Range

Department of Mathematics

State University of New York at Albany

range@albany.edu

\end{document}